\documentclass[12pt,twoside,final,psamsfonts]{amsart}

\usepackage[psamsfonts]{amssymb}
\usepackage{times,a4wide}

\theoremstyle{plain}
\newtheorem*{theorem*}{Theorem}

\newtheorem{theorem}{Theorem}[section]
\newtheorem{proposition}[theorem]{Proposition}
\newtheorem*{proposition*}{Proposition}
\newtheorem{corollary}[theorem]{Corollary}
\newtheorem*{corollary*}{Corollary}
\newtheorem{lemma}[theorem]{Lemma}
\newtheorem*{lemma*}{Lemma}

\theoremstyle{definition}
\newtheorem{remark}[theorem]{Remark}
\newtheorem*{remark*}{Remark}
\newtheorem{example}[theorem]{Example}

\theoremstyle{definition}

\newtheorem*{definition*}{Definition}
\newcommand{\nc}{\newcommand}

\newcommand{\N}{{\mathbb N}}
\newcommand{\D}{{\mathbb D}}
\newcommand{\R}{{\mathbb R}}

\newcommand{\T}{{\partial \mathbb D}}

\renewcommand{\a}{{\alpha}}
\newcommand{\be}{{\beta}}

\newcommand{\Int}{\operatorname{Int}}
\newcommand{\Hol}{\operatorname{Hol}}
\newcommand{\Har}{\operatorname{Har}}
\newcommand{\Arg}{\operatorname{Arg}}
\newcommand{\Cone}{\operatorname{Cone}}
\newcommand{\Conv}{\operatorname{Conv}}

\newcommand{\lra}{\longrightarrow}

\newcommand{\eps}{\varepsilon}
\newcommand{\vp}{\varphi}
\nc{\bea}{\begin{eqnarray}}
\nc{\eea}{\end{eqnarray}}
\nc{\beqa}{\begin{eqnarray*}}
\nc{\eeqa}{\end{eqnarray*}}
\nc{\Hi}{H^{\infty}}
\nc{\loi}{\ell^{\infty}}
\nc{\NL}{N^+\vert \Lambda}
\nc{\hf}{{\mathcal H}_{\phi}}
\nc{\liL}{\lambda\in\Lambda}
\nc{\nn}{\nonumber}

\newenvironment{proof*}{\vskip 2mm\noindent {}}{$\blacksquare$ \vskip 2mm}
\numberwithin{equation}{section}

\newcommand{\Dd}{{\mathbb{D}}}

\renewcommand{\a}{\alpha}

\newcommand{\la}{\lambda}

\newcommand{\eit}{e^{i\theta}}

\renewcommand{\Re}{\mbox{Re}}
\renewcommand{\Im}{\mbox{Im}}

\renewcommand{\qedsymbol}{$\blacksquare$}

\title[Interpolation in the Nevanlinna class and harmonic majorants]
{Interpolation in the Nevanlinna class and harmonic majorants}

\author[A. Hartmann, X. Massaneda, A. Nicolau, P. Thomas]{Andreas Hartmann,
Xavier Massaneda, Artur Nicolau, \& Pascal Thomas}

\address{Laboratoire de Math\'ematiques Pures de Bordeaux,
Universit\'e Bordeaux I, 351 cours de la Lib\'eration,
33405 Talence, France}

\address{Departament de Matem\`atica Aplicada i An\`alisi,
Universitat  de Barcelona, Gran Via 585, 08071-Bar\-ce\-lo\-na, Spain}

\address{Departament de Matem\`atiques, Facultat de Ci\`encies,
Universitat Aut\`onoma de Barcelona, 08193-Bellaterra, Spain}

\address{Laboratoire de Math\'ematiques Emile Picard, UMR CNRS 5580,
Universit\'e Paul Sabatier, 118 route de Narbonne, 31062 TOULOUSE
CEDEX, France.}

\thanks{All authors supported by the PICS program no 1019 of Generalitat de
Catalunya and CNRS. First and third author also supported  by European
Commission Research Training Network HPRN-CT-2000-00116. Second author also
supported by the DGICYT grant BFM2002-04072-C02-01 and the CIRIT grant
2001-SGR00172. Third author also supported by DGICYT grant BFM2002-00571 and
the CIRIT grant 2001-SGR00431}

\email{hartmann@math.u-bordeaux.fr, xavier@mat.ub.es, artur@mat.uab.es,
pthomas@cict.fr}

\date{\today}



\begin{document}

\begin{abstract} We consider a free interpolation problem in Nevanlinna and
Smirnov classes and find a characterization of the corresponding interpolating
sequences in terms of the existence of harmonic majorants of certain functions.
We also consider the related problem of characterizing positive functions in
the disc having a harmonic majorant. An answer is given in terms of a dual
relation which involves positive measures in the disc with bounded Poisson
balayage. We deduce necessary and sufficient geometric conditions, both
expressed in terms of certain maximal functions.
\end{abstract}

\maketitle

\tableofcontents

\section{Introduction and statement of results}
\label{intro}

\subsection{Interpolating sequences for the Nevanlinna Class} Let $\Lambda$ be
a discrete sequence of points in the unit disk $\D$. For a space of holomorphic
functions $X$, the interpolation problem consists in describing the trace of
$X$ on $\Lambda$, i.e.\ the set of restrictions $X\vert \Lambda$, regarded as a
sequence space. One approach is to fix a target space $l$ and look  for
conditions so that $X\vert \Lambda=l$. An alternative approach, known as free
interpolation, is  to require that $X\vert \Lambda$ be {\it ideal}, i.e. stable
under multiplication by $\ell^\infty$. See \cite[Section  C.3.1 (Volume
2)]{Nik02}, in particular, Theorem C.3.1.4, for functional analytic
motivations. This approach is natural for those spaces that are stable under
multiplication by $\Hi$, the space of bounded holomorphic functions on $\D$.
For Hardy and Bergman spaces both definitions turn out to be equivalent, with
the usual choice of $l$ as an $\ell^p$ space with an appropriate weight (see
\cite{ShHSh}, \cite{Se93}).

The situation changes for the non-Banach classes we have in
mind, namely the \emph{Nevanlinna class}
\[
N=\bigl\{f\in \Hol(\D):\lim_{r \to 1}\frac{1}{2\pi}\int_{0}^{2\pi}
    \log^+|f(re^{i\theta})|\;d\theta<\infty\bigr\}
\]
and the related \emph{Smirnov class }
\[
    N^+=\bigl\{f\in N:\lim_{r \to 1}\frac{1}{2\pi}\int_{0}^{2\pi}
    \log^+|f(re^{i\theta})|\;d\theta=\frac{1}{2\pi}\int_{0}^{2\pi}
    \log^+|f(e^{i\theta})|\;d\theta\bigr\}.
\]

We briefly discuss the known results. Naftalevi\v c \cite{Na56} described the
sequences $\Lambda$ for which the trace $N\vert \Lambda$ coincides with the
sequence space $l_{\text{Na}}=\{(a_\lambda)_\lambda: \sup_{\lambda}
(1-|\lambda|)\log^+|a_\lambda|<\infty\}$ (we state the precise result after
Proposition \ref{thmnafta}). The choice of $l_{\text{Na}}$ is motivated by the
fact that $\sup_z (1-|z|)\log^+ |f(z)|<\infty$ for $f\in N$, and this growth is
attained. Unfortunately, the growth condition imposed in $l_{\text{Na}}$ 
forces the sequences to be confined in a finite union of Stolz angles. 
Consequently a big class of Carleson sequences (i.e.\ sequences such that
$\Hi|\Lambda=\ell^{\infty}$), namely those containing a subsequence tending
tangentially to the boundary,  cannot be interpolating  in the sense of
Naftalevi\v c. This does not seem  natural, for $\Hi$ is in the multiplier
space of $N$.  In a sense, the target space $l_{\text{Na}}$ is ``too big".
Further comments on Naftalevi\v c's result can be found in \cite{HM2} and
below, after Proposition \ref{thmnafta}.

For the Smirnov class, Yanagihara \cite{yana2} proved that in order that
$N^+\vert \Lambda$  contains the space
$l_{\text{Ya}}=\{(a_\lambda)_\lambda:\sum_\lambda (1-|\lambda|)
\log^+|a_\lambda|<\infty\}$,  it is sufficient that $\Lambda$ is a Carleson
sequence. However there are Carleson sequences such that $N^+\vert \Lambda$
does not embed into  $l_{\text{Ya}}$ \cite[Theorem 3]{yana2} : the target space
$l_{\text{Ya}}$ is ``too small".

We now turn to the definition of free interpolation.

\begin{definition*}\label{ideal} A sequence space $l$ is called {\it ideal} if
$\ell^\infty l\subset l$, i.e.\ whenever $(a_n)_n\in l$ and  $ (\omega_n)_n\in
\ell^\infty$,  then also $(\omega_n a_n)_n\in l$.
\end{definition*}

\begin{definition*}\label{freeint} Let $X$ be a space of holomorphic functions
in $\D$. A sequence $\Lambda\subset \D$ is called {\it free interpolating for
$X$} if $X\vert \Lambda$ is ideal. We denote $\Lambda\in \Int X$.
\end{definition*}

\begin{remark}
\label{l-infinit}
For any function algebra $X$ containing the constants,
$X\vert\Lambda$ is ideal if and only if
\beqa
    \loi\subset X\vert\Lambda.
\eeqa

The inclusion is obviously necessary. In order to see that it is sufficient
notice that, by assumption, for any $(\omega_\lambda)_\lambda\in \loi$ there
exists $g\in X$ such that $g(\lambda)=\omega_\lambda$. Thus, if
$(f(\lambda))_\lambda\in X\vert \Lambda$, the sequence of values
$(\omega_\lambda f(\lambda))_\lambda$ can be interpolated by  $fg\in X$.
\end{remark}

It is then clear that $\Int N^+\subset\Int N$.

Since the Nevanlinna and Smirnov classes contain the constants, free
interpolation for these classes entails the existence of a nonzero function
$f\in N$ vanishing on $\Lambda$, hence the Blaschke condition $\sum_{\liL}
(1-|\lambda|)<\infty$ is necessary and will be assumed throughout this paper.

Given the Blaschke product $B_\Lambda=\prod_{\lambda\in\Lambda}
b_{\lambda}$ with zero-sequence $\Lambda$,
denote $B_\lambda=B_{\Lambda\setminus \{\lambda\}}=B_\Lambda/b_\lambda$.
Here $b_\lambda =(|\lambda|/\lambda)
(\lambda-z)(1-\overline{\lambda}z)^{-1}$. Define then
\[
\varphi_\Lambda (z) :=
\begin{cases}
\log |B_\lambda(\lambda)|^{-1}\quad&\textrm{if $z=\lambda \in \Lambda$}\\
\ 0 \quad&\textrm{if $z \notin \Lambda$}
\end{cases}
\]

\begin{definition*}
We say that a Borel measurable function $\varphi$ defined on the unit disc
admits a \emph{positive harmonic majorant} if and only if there exists a
positive harmonic function $h$ on the unit disc such that $h(z)\ge \varphi(z)$
for any $z \in \mathbb D$.
\end{definition*}

Let $\Har(\D)$ denote the space of harmonic functions in $\D$ and
$\Har_+(\D)$ the subspace of its positive functions. Consider also
the Poisson kernel in $\D$:
\[ 
P(z,\zeta)= P_z (\zeta)=
\textrm{Re}\left(\frac{\zeta+z}{\zeta-z}\right)=\frac{1-|z|^2}{|\zeta-z|^2}.
\]

Our characterization of interpolating sequences for the Nevanlinna class is as
follows. Note that the existence of a harmonic majorant occurs at two
junctures: first, to decide which sequences of points are free interpolating,
second, to identify the trace space that arises for those sequences which are
indeed free interpolating.

\begin{theorem}\label{thmCNSN}
Let $\Lambda$ be a sequence in $\D$. The following statements are
equivalent:
\begin{itemize}
\item[(a)] $\Lambda$ is a free interpolating sequence for the Nevanlinna
class $ N$.

\item[(b)] The trace space is given by:
\[
   N|\Lambda=l_{N}:=\{ (a_{\lambda})_{\lambda}: \ \exists\
h\in\Har_+(\D)\
\text{such that $h(\lambda)\ge \log^+|a_{\lambda}|$, $\liL$} \}.
\]

\item[(c)] $\varphi_\Lambda$ admits a harmonic majorant.

\item[(d)] There exists $C>0$ such that
for any sequence of nonnegative numbers $\{c_\lambda\}$,
\[
\sum_{\liL} c_\lambda \varphi_\Lambda(\lambda) =
\sum_{\liL} c_\lambda \log |B_\lambda(\lambda)|^{-1}
\le
C \sup_{\zeta \in \partial \mathbb D} \sum_{\liL}  c_\lambda P_{\lambda}
(\zeta).
\]
\end{itemize}
\end{theorem}

We recall that any positive harmonic function on the unit disk is the Poisson
integral of a positive measure on the unit circle,
\[
h(z) = P[\mu](z) = \int_{\partial \mathbb D} P_z(\zeta) \, d\mu(\zeta).
\]

We will say that a harmonic function is \emph{quasi-bounded} if and only if it
admits an absolutely continuous boundary measure (for the reasons for this
terminology, see \cite[pp. 6--7]{He69}). The analogous result for the Smirnov
class will, as can be expected, involve quasi-bounded harmonic
functions.

\begin{theorem}\label{thmCNS}
Let $\Lambda$ be a sequence in $\D$. The following statements are
equivalent:
\begin{itemize}
\item[(a)] $\Lambda$ is a free
interpolating sequence for the Smirnov class $N^+$.

\item[(b)] The trace space is given by
\[
   N^+|\Lambda=l_{N^+}:=\{(a_{\lambda})_{\lambda} : \ \exists\
h\in\Har_+(\D)\
\text{ quasi-bounded }\ : \
h(\lambda)\ge \log^+|a_{\lambda}|, \ \liL \}.
\]

\item[(c)] $\varphi_\Lambda$ admits a quasi-bounded harmonic
majorant.
\end{itemize}
\end{theorem}

The classical Carleson condition characterizing interpolating sequences for
bounded analytic functions in the unit disc is $\sup_\D \varphi_\Lambda  <
\infty$, hence statements (c) in both results above can be viewed as
Carleson-type conditions.

In view of Theorems~\ref{thmCNSN} and \ref{thmCNS}, it seems natural to ask
whether the measure $\mu$ such that $\varphi_\Lambda\leq P[\mu]$ can be
obtained from $\Lambda$ in a canonical way. We do not have an answer to this
question, but with Propositions~\ref{thmnafta} and \ref{supersep} it is easy to
construct examples that discard natural candidates, such as the (weighted) sum
of Dirac masses $\mu=\sum_\lambda (1-|\lambda|)\delta_\lambda$, or its Poisson
balayage (see definition below) $d\mu=\sum_\lambda (1-|\lambda|)
P_\lambda(\zeta) d\sigma$.

\subsection{Positive harmonic majorants}
The condition in Theorem \ref{thmCNSN} (d)
arises in the solution of a problem of independent interest:

{\bf Problem.} {\em Which functions $\varphi : \mathbb D \longrightarrow
\mathbb R_+$ admit a (positive) harmonic majorant?}

Answers to this problem lead to rather precise theorems about the permissible
decrease of the modulus of bounded holomorphic functions, e.g. Corollary
\ref{decrease} below. See  \cite{Hay}, \cite{LySe97}; \cite{EiEs} also provides
a survey of such results. The existence of harmonic majorants is relevant as
well to the study of zero-sequences for Bergman and related spaces of
holomorphic functions \cite{Lu96}.

An answer to the problem of positive harmonic majorants can be given in dual
terms (see \cite{NiTh} for another, more theoretical, characterization). The
\emph{Poisson balayage} (or swept-out function) of a finite positive measure
$\mu$ in the closed unit disc is defined as
\[
B(\mu) ( \zeta) = \int_{\mathbb D} P_z (\zeta) d \mu (z)
\]
for $\zeta \in \partial \mathbb D$.
We will be interested in the class of measures having bounded
balayage. Recall that Carleson measures are those finite positive
measures whose balayage has bounded mean oscillation \cite[Theorem VI.1.6,
p. 229]{Gar}; this is also an easy consequence of the $H^1$-BMO duality
\cite[Theorem VI.4.4, p. 245]{Gar}.
Hence positive measures with bounded balayage form a subclass of
the usual Carleson measures. It is easy to see
(cf. Section \ref{secFarkas}) that positive
measures with bounded balayage are precisely those which operate
against positive harmonic functions, that is, those measures $\mu$
for which there exist a constant $C=C(\mu)$ such that
\[
\int_{\mathbb D} h(z) d \mu(z) \le C h(0)
\]
for any positive harmonic function in the unit disc $\mathbb D$.

\begin{theorem}
\label{thmequivBB}
Let $\varphi$ be a nonnegative Borel  function on
the unit disk $\mathbb D$. Then there exists a (positive) harmonic
function $h$ such that $h (z) \ge \varphi(z)$ for any $z \in \D$
if and only if there exists a constant $C=C(\varphi)$ such that
for every positive measure $\mu$,
\begin{equation*}
\int_{\D}  \varphi (z)\, d\mu (z) \le C \sup_{\zeta
\in \partial\mathbb D} \int_{\D} P_z( \zeta) \, d\mu(z).
\end{equation*}
\end{theorem}

The necessity of the condition is obvious, while the sufficiency follows
from a convenient version of a classical result in Convex Analysis,
known as Minkowski-Farkas Lemma. The characterization of interpolating
sequences in the Nevanlinna class in dual terms given by condition (d)
in Theorem \ref{thmCNSN} follows from this result.

This can be applied to study the decrease of a non-zero bounded analytic
function in the disc along a given  non-Blaschke sequence.

\begin{corollary}\label{decrease} Let $\Lambda$ be a separated non-Blaschke
sequence and $(\eps_\lambda)_{\liL}$ a sequence of positive values. Then there
exists a non-zero  function $f\in H^\infty(\D)$ with $|f(\lambda)| <
\varepsilon_\lambda$, $\liL$, if and only if  $\Lambda$ is the union of a
Blaschke sequence and a sequence $\Gamma$ for which there exists a universal
constant $C=C(\Gamma)$ such that
\[
\sum_{\gamma \in \Gamma} c_\gamma \log \varepsilon_\gamma^{-1}   \le C
\sup_{\zeta \in \partial \D }
\sum_{\gamma \in \Gamma} c_\gamma P_\gamma  (\zeta)
\]
for any sequence of nonnegative numbers $(c_\gamma)_{\gamma\in\Gamma}$.
\end{corollary}

We do not know what is the corresponding analogue of Theorem~\ref{thmequivBB}
for quasi-bounded harmonic functions (i.e. for the Smirnov class). In proving
the sufficiency of such an analogue, one should deduce that the measure $\mu$
giving $h=P[\mu]$ is absolutely continuous with respect to the Lebesgue
measure. This extra condition on $\mu$ does not seem to follow in an obvious
way from our proof, based on Farkas' Lemma.

For the problem of harmonic majorants it is desirable to obtain
criteria which, although only necessary or sufficient, are more geometric and
easier to check than the duality condition of Theorem~\ref{thmequivBB}.

Recall that the Stolz angle with vertex  $\zeta \in \partial\D$ and
aperture $\alpha$ is defined by
$$
\Gamma_\alpha (\zeta) := \{ z \in \mathbb D : |z-\zeta| \le \alpha
(1-|z|^2) \}.
$$
In our considerations the angle $\alpha$ is of no importance, so
we will write $\Gamma(\zeta)$ for the generic Stolz angle with
aperture $\alpha$. Given a function $f$ from $\D$ to $\mathbb
R_+$, the non-tangential maximal function is defined as
\[
M f (\zeta) := \sup_{\Gamma (\zeta)} f .
\]
Let $\sigma$ be the normalized Lebesgue measure on $\partial\D$.
Consider the weak-$L^1$ space
\[
L_w^{1}(\partial\D)=\{f\ \textrm{measurable}\, :
\sup_{t>0} t \sigma\{\zeta : |f(\zeta)|>t\}<\infty\}
\]
and let
\[
L^{1}_{w,0}(\T)=\{ f\ \textrm{measurable}\,: \lim_{t\to \infty}t
\sigma\{\zeta : |f(\zeta)|>t\}=0\}.
\]

It is well-known that the non-tangential maximal function of the Poisson
transform of a positive finite measure belongs to $L^{1}_{w}$ \cite[Theorem
5.1, p. 28]{Gar}.  A more careful analysis shows that if $\mu$ is  absolutely
continuous, then its Poisson transform is in $L^{1}_{w,0}$. This and some easy
estimates imply the following result.

\begin{proposition}\label{maxfcnharm}
\begin{itemize}

\item[(a)] If $\varphi$ admits a harmonic majorant, then $M\varphi \in
L^1_w(\partial \mathbb D)$.

\item[(b)] If $\varphi$ admits a positive quasi-bounded harmonic majorant, then
$M\varphi \in L^1_{w,0}(\partial \mathbb D)$.

\item[(c)] If $M \varphi\in L^1 (\partial \mathbb D)$, then
the function $\varphi$ admits $P[M\varphi]:=P[M \varphi  d\sigma]$
as a quasi-bounded harmonic majorant.
\end{itemize}
\end{proposition}

As far as necessary conditions are concerned, there is a way to improve
the previous result by using another maximal function.

\begin{definition*} Given a nonnegative function  $f$ defined on the
unit circle, its \emph{Hardy-Littlewood maximal function} is given by
$$
f^*(x) := \sup_{I \ni x} \frac1{\sigma(I)} \int_I f.
$$
\end{definition*}

For $\varphi$ a nonnegative function on the disc,
define
\[
\varphi^H (\zeta) := \sup_{z \in \D}  \varphi(z) \chi_{I_z}^*(\zeta),
\]
where $I_z$ is the ``Privalov shadow" interval
\begin{equation}
\label{shadow}
I_z := \{ \zeta \in \partial \D : z \in \Gamma(\zeta) \}.
\end{equation}

\begin{proposition}
\label{strongneccond}
\begin{itemize}
\item[(a)] If $\varphi$ admits a harmonic majorant, then
$\varphi^H \in L^1_w(\partial \D) $.
\item[(b)] If $\varphi$ admits a quasi-bounded harmonic
majorant, then $\varphi^H \in L^1_{w,0}(\partial \D) $.
\end{itemize}
\end{proposition}

We will give some examples in Proposition~\ref{notsuffstronger} that show that
this is indeed stronger than the necessary condition given in the first part of
Proposition \ref{maxfcnharm}, but still falls short of giving a sufficient 
condition for the existence of a harmonic majorant.

\subsection{Geometric criteria for interpolation} We would like to obtain some
geometric implications of the analytic conditions given in Theorems
\ref{thmCNSN} and \ref{thmCNS}. To begin with, we would like to state the maybe
surprising result that separated sequences (with respect to the hyperbolic
distance) are interpolating for the Smirnov class (and hence the Nevanlinna
class). Recall that a sequence $\Lambda$ is called \emph{separated} if
$\delta(\Lambda):=\inf\limits_{\lambda\neq \lambda'} \rho(\lambda,\lambda')
>0$, where
$$
\rho(z,w):=|b_z(w)|= \left| \frac{z-w}{1-z\overline w} \right|,
$$
is the pseudo-hyperbolic distance.

For such sequences, the values $\log |B_\lambda(\lambda)|^{-1}$ can always be
majorized by the values at $\lambda \in \Lambda$ of the Poisson integral of an
integrable function (see Proposition \ref{propsep}), thus the following
corollary is immediate from Theorem~\ref{thmCNS}.

\begin{corollary}\label{separated}
Let $\Lambda$ be a separated Blaschke sequence. Then
$\Lambda\in\Int N^+$ (hence $\Lambda\in \Int N$).
\end{corollary}

More precise conditions can be deduced from Propositions \ref{maxfcnharm},
\ref{strongneccond} and (c) in Theorems \ref{thmCNSN} and \ref{thmCNS}.

\begin{corollary}\label{maxfct}
Let $\Lambda $ be a sequence in $\D$.

\begin{itemize}
\item[(a)] If $\Lambda\in\Int N$ then $\varphi_\Lambda^H\in L^1_{w}(\T)$.
If $\Lambda\in\Int N^+$ then $\varphi_\Lambda^H\in L^1_{w,0}(\T)$.

\item[(b)] If $M\varphi_\Lambda\in L^1(\T)$ then $\Lambda\in\Int N^+$
(and hence $\Lambda\in \Int N$).
\end{itemize}
\end{corollary}

Notice that the necessary conditions obtained by replacing  $
\varphi_\Lambda^H$ by $M\varphi_\Lambda$ in (a) also hold. This in an
immediate consequence of the estimate $\varphi_\Lambda^H\geq
M\varphi_\Lambda$.

This result implies the following Carleson-type conditions.

\begin{corollary}
\label{CNgeom}
({\rm a}) If $\Lambda\in\Int N^+$, then
\bea\label{CN}
    \lim_{|\lambda|\to 1} (1-|\lambda|)\log  |B_{\lambda}(\lambda)|^{-1}=0.
\eea

({\rm b}) If $\Lambda\in\Int N$, then \bea\label{CNN}
    \sup_{\lambda\in\Lambda} (1-|\lambda|)
    \log|B_{\lambda}(\lambda)|^{-1}<\infty.
\eea

({\rm c}) If $\Lambda$ is Blaschke and
\bea\label{CS}
    \sum_{\lambda\in\Lambda} (1-|\lambda|)
    \log|B_{\lambda}(\lambda)|^{-1}<\infty ,
\eea then $\Lambda\in \Int N^+$ (and so also $\Lambda\in \Int N$).

\end{corollary}

Condition~\eqref{CN} already appeared in \cite[Theorem 1]{yana2} as a necessary
condition for the sequence space $l_{\text{Ya}}$ (as defined in the beginning
of Section 1.1) to be included in the trace of $N^+$.

In some situations the conditions above are indeed a characterization of
interpolating sequences. For instance,  the weak $L^1$-condition characterizes
interpolating sequences lying on a radius, while for sequences approaching the
unit circle very tangentially the characterization is given by the strong
$L^1$-condition. This is collected in the next results.

\begin{proposition}\label{thmnafta}
Assume that $\Lambda\subset\D$ lies in a finite union of Stolz
angles.
\begin{itemize}
\item[(a)] $\Lambda\in\Int N^+$ if and only if \eqref{CN} holds.

\item[(b)] $\Lambda\in\Int N$ if and only if \eqref{CNN} holds.

\end{itemize}
\end{proposition}

It should be mentioned that (b) can also be derived from Naftalevi\v c's result
\cite[Theorem 3]{Na56}. On the other hand, his full characterization of the
sequences such that $N|\Lambda = l_{\text{Na}}$ can also be deduced from
Theorem \ref{thmCNSN}.

\begin{corollary*}[Naftalevi\v c, 1956] $N\vert \Lambda = l_{\text{Na}}$ if and
only if $\Lambda$ is contained in a finite union of Stolz angles and
\eqref{CNN} holds.
\end{corollary*}

Let us consider
the other geometric extreme, sequences which in particular only approach the
circle in a tangential fashion. Write
\begin{equation}\label{eqnp7}
\mu_\Lambda := \sum_\la (1-|\la|) \delta_\la,
\end{equation}
where $\delta_\la$ stands for the Dirac measure at $\la$.

\begin{proposition}
\label{supersep}
If $\mu_\Lambda$ has bounded balayage, then $\Lambda\in\Int N$ if and only if
$\Lambda\in\Int N^+$, and if and only if \eqref{CS} holds.
\end{proposition}

Note that the condition that $\mu_\Lambda$ has bounded balayage implies in
particular that $\Lambda$ approaches the circle tangentially. In Section
\ref{GeomCond}, we will see more concrete conditions of geometric separation
which are sufficient to imply that $\mu_\Lambda$ has bounded balayage
(Proposition \ref{thmnafta2}).

When $\mu_\Lambda$ has bounded balayage, 
the trace space will embed into Yanagihara's target space. 
More precisely, the following result holds.

\begin{proposition}
\label{embedyana}
The following are equivalent:
\begin{itemize}
\item[(a)] $N|\Lambda \subset l_{\text{Ya}}$,

\item[(b)] $N^+|\Lambda \subset l_{\text{Ya}}$,

\item[(c)] $\mu_\Lambda$ has
bounded balayage, i.e. 
$\sup_{\zeta\in\partial\D}\sum_\la (1-|\la|) P_\la(\zeta) <\infty$.
\end{itemize}
\end{proposition}

Yanagihara considered the sequences such that $N^+|\Lambda \supset
l_{\text{Ya}}$. These are automatically in $\Int N^+$, since for any Blachke
sequence $l_{\text{Ya}} \supset \ell^\infty$. Conversely, Lemma
\ref{L1suffcond} (see Section \ref{GeomCond}) implies that $l_{\text{Ya}}
\subset l_{N^+}$, thus if $\Lambda \in \Int N^+$, then  by Theorem
\ref{thmCNS}(b) $N^+|\Lambda \supset l_{\text{Ya}}$. Therefore Theorem
\ref{thmCNS} characterizes in particular the sequences studied by Yanagihara.

Altogether, free interpolation for the Nevanlinna and Smirnov classes can be
described in terms of the intermediate target spaces $l_N$ and $l_{N^+}$.
Notice first that always
$N^+|\Lambda \subset l_{N^+}$ and $N|\Lambda \subset l_{N}$ (this is proved
at the beginning of Section \ref{trace}). So,
$\Lambda\in\Int N^+$ if and only if $N^+|\Lambda \supset l_{N^+}$,
and $\Lambda\in\Int N$ if and only if $N|\Lambda
\supset l_{N}$. 
Observe also that $l_{\text{Ya}} \subset
l_{N^+} \subset l_{N} \subset l_{\text{Na}} $.

The paper is organized as follows. The next section is devoted to collecting
some basic results on functions in the Nevanlinna class. In Section
\ref{condsuff} we prove the sufficiency for interpolation of the conditions (c)
of Theorems~\ref{thmCNS} and \ref{thmCNSN}. We essentially use a result by
Garnett allowing interpolation by $\Hi$ functions on sequences which are denser
than Carleson sequences, under some decrease assumptions on the interpolated
values. In Section~\ref{NECESSITY} we study the necessity of these conditions.
We first observe that in the product  $B_{\lambda} (\lambda)$ appearing in
Theorem ~\ref{thmCNSN}, only the factors $b_{\lambda}(\lambda')$ with
$\lambda'$ close to $\lambda$ are relevant. Then we split the sequence into
four pieces, thereby reducing the interpolation problem, in a way, to that on
separated sequences. The trace space characterization will be discussed in
Section ~\ref{trace}. In Section ~\ref{secFarkas} we consider measures with
bounded balayage, show that they operate against positive harmonic functions
and prove Theorem ~\ref{thmequivBB}. In Section \ref{weakcond}, we prove
Proposition \ref{strongneccond}, and provide examples to show that the
sufficient condition is not necessary, and the necessary condition not
sufficient. Section  ~\ref{GeomCond} is devoted to the proofs of
Corollary~\ref{CNgeom}, Propositions~\ref{thmnafta}, \ref{embedyana}, and
\ref{supersep}, as well as the deduction of Naftalevi\v c's result from Theorem
\ref{thmCNSN}. Also, we give examples of measures with bounded balayage. In the
final section, we exploit the reasoning of Section ~\ref{condsuff} to construct
non-Carleson interpolating sequences for ``big'' Hardy-Orlicz classes.

{\bf Acknowledgements.} The authors wish to express special thanks to
Jean-Baptiste Hiriart-Urruty for introducing them to Farkas' Lemma, and to
Alexander Borichev for Theorem \ref{discrete}.

\section{Preliminaries}\label{prelim}

We next recall some standard facts about the structure of the Nevanlinna and
Smirnov classes (general refe\-rences are e.g. \cite{Gar}, \cite{Nik02} or
\cite{RR}).

A function $f$ is called \emph{outer} if it can be written in the form
\[
f(z)=C \exp \left\{\int_\T \frac{\zeta+z}{\zeta-z}
\log v(\zeta) d\sigma(\zeta) \right\},
\]
where $|C|=1$, $v>0$  a.e.\ on $\T$ and $\log v \in L^1(\T)$. Such a
function is
the quotient $f=f_1/f_2$ of two bounded outer functions  $f_1,f_2\in \Hi$
with  $\|f_i\|_{\infty}\le 1$, $i=1,2$. In particular, the weight $v$
is given by
the boundary values of $|f_1/f_2|$. Setting $w=\log v$,
we have
\[
\log |f(z)|=P[w](z) =\int_\T P_z(\zeta) w(\zeta)d\sigma(\zeta).
\]
This formula allows us to freely switch between assertions about outer
functions $f$ and the associated measures $w d\sigma$.

Another important family in this context are  \emph{inner} functions: $I\in
\Hi$ such that $|I|=1$ almost everywhere on $\T$. Any inner function $I$ can be
factorized into a Blaschke product $B_\Lambda$ carrying the zeros
$\Lambda=\{\lambda_n\}_n$  of $I$, and a singular inner function
$S$ defined by
\beqa
    S(z)=\exp\left\{-\int_\T \frac{\zeta+z}{\zeta-z}\,d\mu(\zeta)\right\},
\eeqa
for some positive Borel measure $\mu$ singular with respect to Lebesgue
measure.

According to the Riesz-Smirnov factorization, any
function $f\in N^+$ is represented as
\[
    f=\alpha \frac{B S f_1}{f_2},
\]
where $f_1, f_2$ are outer with $\|f_1\|_{\infty},\|f_2\|_{\infty}\le 1$,
$S$ is singular inner, $B$ is a Blaschke product and $|\alpha|=1$.
Similarly,
functions $f\in N$ are represented as
\[
    f=\alpha \frac{B S_1 f_1}{S_2 f_2},
\]
with $f_i$ outer, $\|f_i\|_{\infty}\le 1$, $S_i$ singular inner,
$B$ is a Blaschke product and $|\alpha|=1$.

In view of the Riesz-Smirnov factorization described above,  the essential
difference bet\-ween Nevanlinna and Smirnov functions is the extra singular
factor appearing in the denominator in the Nevanlinna case. This is reflected
in the corresponding result for free interpolation in $N$ by the fact that
$\varphi_\Lambda$ is bounded by a harmonic function, not necessarily
quasi-bounded.

\section{From harmonic majorants to interpolation}\label{condsuff}

For a given Blaschke sequence $\Lambda\subset \D$ set
$\delta_\lambda=|B_{\lambda}(\lambda)|$. The key result in the proof of the
sufficient condition is the following theorem by Garnett \cite{Gar77},
that we
cite for our purpose in a slightly weaker form (see also \cite{Nik02} as a
general source, in particular C.3.3.3(g) (Volume 2) for more results of
this
kind).

\begin{theorem*}
Let $\vp:[0,\infty)\longrightarrow [0,\infty)$ be a decreasing
function such that
$\int_0^{\infty}\vp(t)\,dt<\infty$. If a sequence
$a=(a_\lambda)_\lambda$ satisfies
\beqa
|a_\lambda|\le \delta_\lambda\vp (\log\frac{e}{\delta_\lambda}),\quad \liL,
\eeqa
then there exists a function $f\in \Hi$ such that $f|\Lambda=a$.
\end{theorem*}

Observe that according to our former notation we have
$\log (e/\delta_{\lambda})=1+\varphi_{\Lambda}(\lambda)$.

As we have already noted in Remark~\ref{l-infinit}, in order to have free
interpolation in the Nevanlinna and Smirnov classes, it is sufficient that
$\loi\subset N|\Lambda$ and $\loi\subset \NL$ respectively.  Our aim will
be
to accommodate the decrease given in Garnett's result by an appropriate
function in $N$ or $N^+$. This is the crucial step in the proof of  the
sufficient conditions of Theorems~\ref{thmCNS} and \ref{thmCNSN}, and it
occupies its main part.

\begin{proof}[Proof of sufficiency of ~\ref{thmCNS} (c) and ~\ref{thmCNSN} (c)]
The proof will be presented for the more difficult case of the Nevanlinna
class. So, assume that $h\in\Har_+(\D)$ majorizes $\varphi_\Lambda$. Then $h$
is the Poisson integral of a positive measure $\mu$ on the circle and the
function
\begin{equation}
\label{outerrepr}
    g(z)= \int_\T\frac{\zeta+z}{\zeta-z}d\mu(\zeta)
\end{equation}
has positive real part in the disc. By Smirnov's theorem,
$g$ is an outer function in some $H^p$, $p<1$, and therefore in
$N^+$ (see \cite{Nik02}, in particular A.4.2.3 (Volume 1)). Also
$\exp(g)$ is in the Nevanlinna class. By assumption we have $\log
(1/\delta_{\lambda})\le \textrm{Re}\; g(\lambda)$,
$\lambda\in\Lambda$.

Take now $\vp(t)=(1+t)^{-2}$, which obviously satisfies the
hypothesis of Garnett's theorem, and set $H=(2+g)^2$, which is
still outer in $N^+$. We have the estimate
\beqa
    |H(\lambda)| = |2+g(\lambda)|^2
     \ge (2+\textrm{Re}\; g(\lambda))^2
     \ge (1+\log \frac e{\delta_\lambda})^2
    =\frac{1}{\vp(\log(e/\delta_\lambda))},
\eeqa
hence the sequence $(\gamma_\lambda)_\lambda$ defined by
\beqa
    \gamma_\lambda=\frac{1}{H(\lambda)\vp(\log(e/\delta_\lambda))},
    \quad \lambda\in\Lambda,
\eeqa
is bounded by $1$.

In order to interpolate $\omega=(\omega_\lambda)_\lambda\in\loi$ by a
function in $N$, split
\beqa
    \omega_\lambda=
    \big(\frac{\omega_\lambda \gamma_\lambda\exp
(-g(\lambda))}{\delta_\lambda}
    \delta_\lambda\vp(\log\frac  e{\delta_\lambda})\big) \cdot
     \frac{H(\lambda)}{\exp( -g(\lambda))}.
\eeqa

Since by hypothesis $(\omega_\lambda\gamma_\lambda \exp(-g(\lambda))
/\delta_\lambda)_\lambda$ is bounded, we
can apply Garnett's result to interpolate the sequence
\beqa
    a_\lambda=
    \frac{\omega_\lambda\gamma_\lambda \exp(-g(\lambda))}{\delta_\lambda}
    \delta_\lambda
    \vp(\log\frac e{\delta_\lambda}),\quad \lambda\in \Lambda,
\eeqa by a function $f\in\Hi$. Now $F=fH \exp(g)$ is a function in
$N$ with $F|\Lambda=\omega$.

The proof for the Smirnov case is obtained by observing that if
the measure $\mu$ is absolutely continuous, then the function $\exp(g)$
is in the Smirnov class and so is the interpolating function $F$.
\end{proof}

\section{From interpolation to harmonic majorants}
\label{NECESSITY}

We first show that in order to construct the appropriate function
estimating $\log|B_{\lambda}(\lambda)|^{-1}$ we only need to consider
the  factors of $B_\lambda$ given by points $\lambda'\in\Lambda$ which
are close to $\lambda$. This is in  accordance with the results for some
related spaces of functions \cite[Theorem 1]{HM2}, and it obviously
implies Corollary~\ref{separated}.

\begin{proposition}
\label{propsep}
Let $\Lambda$ be a Blaschke sequence.  For any $\delta \in (0,1)$, there
exists a quasi-bounded positive harmonic function $h=P[w]$,
$w \in L^1(\partial \D)$, such that
\beqa
    -\log \prod_{{\lambda:\rho(\lambda,z)\ge \delta}} |b_{\lambda}(z)|
    \leq h(z) ,\quad
    z \in \D ,
\eeqa
and therefore an outer function $G\in N^+$, where $G=\exp(-g)$
and $g$ is given by \eqref{outerrepr} with
$d\mu=w d\sigma$, such that
\beqa
    \prod_{{\lambda:\rho(\lambda,z)\ge \delta}} |b_{\lambda}(z)|
    \ge |G(z)|,\quad
    z \in \D.
\eeqa
\end{proposition}


\begin{proof}
We shall use the intervals $I_z$ introduced in
\eqref{shadow}.
In \cite[p. 124, lines 3 to 17]{NPT}, it is proved that the
function $w$  given by
$$
w(\zeta) = c_0 \sum_{\lambda \in \Lambda} \chi_{I_\lambda} (\zeta),
$$
where $c_0$ is an appropriate positive
constant, is suitable.
At this juncture, the separation hypothesis
made in \cite[Lemma 4]{NPT} is no longer used.
\end{proof}

\begin{proof}[Proof of the necessity of \ref{thmCNS} (c) and \ref{thmCNSN}
(c)]
We will use a dyadic partition of the disk:
for any $n$ in $\mathbb N$, let
\begin{equation}
\label{dyadarcs}
I_{n,k} := \{ \eit : \theta \in [{2\pi}k 2^{-n},{2\pi}(k+1)
2^{-n}) \}, \ 0 \le k < 2^n .
\end{equation}
and the associated Whitney partition in ``dyadic squares":
\begin{equation}
\label{dyadsquares}
Q_{n,k} := \{ r \eit : \eit \in I_{n,k} , 1- 2^{-n} \le r < 1-
2^{-n-1} \}.
\end{equation}
Observe that the hyperbolic diameter of each Whitney square
$Q_{n,k}$ is bounded between two absolute constants. 

We split the sequence into four pieces: $\Lambda=
\bigcup_{i=1}^4 \Lambda_i$ such that each piece $\Lambda_i$ lies in a union
of dyadic squares
that are uniformly separated from each other.
More precisely, set
\beqa
    \Lambda_1=\Lambda\cap Q^{(1)},
\eeqa
where the family $Q^{(1)}$ is given by $\{Q_{2n,2k}\}_{n,k}$
(for the remaining three sequences we respectively choose
$\{Q_{2n ,2k+1}\}_{n,k}$,
$\{Q_{2n+1,2k}\}_{n,k}$ and
$\{Q_{2n+1,2k+1}\}_{n,k}$).
In order to avoid technical difficulties we count only those
$Q$ containing points of $\Lambda$ (in case $\Lambda \cap Q$
is empty there is nothing to prove).
In what follows we will argue on one sequence, say $\Lambda_1$.
The arguments are the same for the other sequences.

Our first observation is that, by construction, for $Q, L \in Q^{(1)}$,
$Q \neq L$,
\beqa
    \rho(Q,L):=\inf_{z\in  Q, w\in L}\rho(z,w)
    \ge\delta>0,
\eeqa
for some fixed $\delta$. In what follows, the letters $j$, $k$... will
stand for indices in $\mathbb N^2$ of the form $(n,l), 0 \le l <2^n$.
The closed rectangles
$\overline{Q_j}$ are compact in
$\D$ so that $\Lambda_1\cap Q_j$ can only
contain a finite number of points (they contain at least one point,
by assumption).
Therefore
\beqa
    0<m_j :=\min_{\lambda\in \Lambda_1\cap Q_j} |B_{\lambda}(\lambda)|
\eeqa
(note that we consider the entire Blaschke product $B_\lambda$
associated with $\Lambda\setminus\{\lambda\}$). Take
$\lambda_j^{1}\in Q_j$ such that
$m_j=|B_{\lambda_j^{1}}(\lambda_j^{1})|$.

Assume now that $\Lambda\in \Int N$.
Since $\loi\subset N\vert\Lambda$, there exists a function $f_1\in N$
such that
\beqa
    f_1(\lambda)=
    \begin{cases}
1\quad &\text{if}\quad \lambda\in \{\lambda_j^{1}\}_j\\
0\quad &\text{if}\quad \lambda\notin \{\lambda_j^{1}\}_j.
\end{cases}
\eeqa
By the Riesz-Smirnov factorization we have
\begin{equation}
\label{factorizn}
    f_1=B_{\Lambda\setminus \{\lambda_j^{1}\}_j}\frac{h_1}{h_2 T_2},
\end{equation}
where $T_2$ is singular inner, $h_1$ is some function in $\Hi$ and $h_2$ is
outer in $\Hi$.
Again, we can assume
$\|h_i\|_{\infty}\le 1$, $i=1,2$. Hence
\beqa
    1 = |f_1(\lambda_k^{1})|
     \le |B_{\Lambda\setminus \{\lambda_j^{1}\}_j}(\lambda_k^{1})|\cdot
      \frac{1}{|h_2(\lambda_k^{1}) T_2(\lambda_k^{1})|},
\eeqa
and
\beqa
    |B_{\Lambda\setminus \{\lambda_j^{1}\}_j}(\lambda_k^{1})|
     \ge {|h_2(\lambda_k^{1}) T_2(\lambda_k^{1})|},\quad  k\in\N.
\eeqa

Since $h_2 T_2$ does not vanish and is bounded above by 1, the function
$\log |h_2 T_2|$ is a negative harmonic function. By Harnack's
inequality, there exists an absolute constant $c\ge 1$ such that
\beqa
    \frac{1}{c} |\log |h_2(\lambda_k^{1}) T_2(\lambda_k^{1})||
    \le |\log |h_2(z) T_2(z)|| \le c |\log|h_2(\lambda_k^{1})
T_2(\lambda_k^{1})||,
    \quad z\in Q_k,
\eeqa
hence
\beqa
|h_2(\lambda_k^{1})T_2(\lambda_k^{1})|^{c}\le |h_2(z) T_2(z)|
    \le |h_2(\lambda_k^{1}) T_2(\lambda_k^{1})|^{1/c},\quad z\in Q_k.
\eeqa
This yields
\bea\label{estimate4}
    |(h_2 T_2)^c(\lambda')|\le |(h_2 T_2)(\lambda_k^{1})|\le
    |B_{\Lambda\setminus \{\lambda_j^{1}\}_j}(\lambda_k^{1})|
\eea
for every $\lambda'\in \Lambda_1\cap Q_k$.

Let us now exploit Proposition \ref{propsep}. By construction,
the sequence $\{\lambda_j^{1}\}_j\subset\Lambda_1$ is
separated. Therefore, there
exists an outer function $G_1$ in the Smirnov class such that
\beqa
    |B_{\{\lambda_j^{1}\}_j\setminus \{\lambda_k^{1}\}}(\lambda_k^{1})|
    \ge |G_1(\lambda_k^{1})|,\quad k\in\N.
\eeqa
Again, $G_1$ is a quotient of two bounded outer functions and
we can suppose that $G_1$ is outer in $\Hi$ with
$\|G_1\|_{\infty}\le 1$. Also, we can use Harnack's inequality as above
to get
\beqa
    |G_1(\lambda_k^{1})|\ge |G_1^c(\lambda')|
\eeqa
for every $\lambda'\in \Lambda_1\cup Q_k$.
This together with \eqref{estimate4} and our definition
of $\lambda_k^{1}$ give
\beqa
    |B_{\Lambda\setminus\{\lambda'\}}(\lambda')|
    &\ge& |B_{\Lambda\setminus\{\lambda_k^{1}\}}(\lambda_k^{1})|
    =|B_{\Lambda\setminus \{\lambda_j^{1}\}_j}(\lambda_k^{1})|
    \cdot |B_{\{\lambda_j^{1}\}_j\setminus \{\lambda_k^{1}\}}
      (\lambda_k^{1})|\\
    &\ge& |(h_2 T_2)^c(\lambda')|\cdot |G_1^c(\lambda')|
\eeqa
for every $\lambda'\in Q_k$ and $Q_k\in Q^{(1)}$. Set  $g_1=(h_2G_1)^c$ and
$S_1=T_2^c$ ; by construction, $g_1$ is outer with $\|g_1\|_{\infty}\le 1$
and
$S_1$ is singular inner.

Construct in a similar way functions $g_i$, $S_i$
for the sequences $\Lambda_i$, $i=2,3,4$ and define the products
$g=\prod_{i=1}^4 g_i$ and $S=\prod_{i=1}^4 S_i$.
Of course $g$ is outer in $\Hi$, and $S$ is singular
inner. So, whenever $\lambda'\in\Lambda$, there exists
$k\in \{1,2,3,4\}$ such that $\lambda'\in\Lambda_k$, and hence
\begin{equation}
\label{blaschkeminor}
    |B_{\lambda}(\lambda)|\ge |g_k(\lambda)S_k(\lambda)|\ge
|g(\lambda)S(\lambda)|.
\end{equation}

Therefore, the positive harmonic function $h=- \log|gS|$ satisfies
$h(\lambda) \ge - \log |B_{\lambda}(\lambda)|$. The proof for
$N^+$ goes along the same lines, except that singular inner factors
do not occur in \eqref{factorizn}, and so will not appear in
\eqref{blaschkeminor} neither.
\end{proof}

\section{The trace spaces}\label{trace}

In this short section we  prove the trace space characterization of free
interpolation given in Theorems \ref{thmCNSN} and \ref{thmCNS}.

In order to see that (b) in each theorem implies free interpolation
it suffices to
observe that  $\ell^{\infty}\subset l_{N^+}\subset  l_N$ and use Remark
\ref{l-infinit}.

For the proof of the converse, we will only consider the situation in
the Nevanlinna class, since the case of the Smirnov class is again
obtained by removing the singular part of the measure and the singular
inner factors.

Assume that $(a_\lambda)_\lambda\in N|\Lambda$ and that $f\in N$ is such
that $f(\lambda)=a_\lambda$, $\lambda\in\Lambda$. Since $f$ can be
written as $f=f_1/(S_2f_2)$, where $f_1\in \Hi$, $\|f_1\|_{\infty}\le
1$, $S_2$ is  singular inner with associated singular measure $\mu_S$,
and $f_2\in \Hi$ is an  outer function with $\|f_2\|_{\infty}\le 1$, we
can define the positive finite measure $\mu=\log (1/|f_2|)\,
d\sigma+d\mu_S$ which obviously satisfies
$P[\mu](\lambda)\geq\log^+|a_\lambda|$, $\lambda\in\Lambda$.

Conversely, suppose that $(a_{\lambda})_{\lambda}$ is
such that there is a positive finite measure $\mu$ with $P[\mu](\lambda)\ge
\log^+|a_{\lambda}|$. The Radon-Nikodym decomposition of $\mu$ is given by
$d\mu=w\,d\sigma+d\mu_S$, where $w\in L^1(\T)$ is positive and
$\mu_S$ is a positive
finite singular measure. Let $S$ be the singular inner function associated
with
$\mu_S$, and let $f$ be the function defined by
\beqa
    f(z)=\exp\left(\int_{\T}\frac{\zeta+z}{\zeta-z}
w(\zeta)\,d\sigma(\zeta)\right),
\quad z\in\D.
\eeqa
By definition, $f$ is outer in $N^+$ and $F=f/S\in N$. Clearly,
$\log^+|a_{\lambda}|\le \log |F(\lambda)|$, thus $|a_{\lambda}|\le
|F(\lambda)|$.
Since $N|\Lambda$ is ideal by assumption,
there exists  $f_0\in N$ interpolating $(a_{\lambda})_{\lambda}$.
\hfill \qedsymbol

\section{Harmonic majorants and measures with bounded
balayage}\label{secFarkas}

Let us start by proving that positive measures with bounded balayage are
precisely those which operate against positive harmonic functions. Define
\begin{equation*}
\mathcal B := \{ \mu \mbox{ positive Borel measures on }\mathbb D
\mbox{ such that } \sup_{\zeta \in \partial \mathbb
D} B(\mu)(\zeta) \le 1  \}.
\end{equation*}

\begin{proposition}
Let $\mu$ be a positive Borel measure on the disk.  Then
$\int_{\D} h d\mu$ is finite for any positive harmonic function
$h$ on the disk if and only if there exists some $c>0$ such that
$\mu$ has balayage uniformly bounded by $c$. Furthermore, the
relevant constants are related:
\[
\sup_{\zeta \in \partial \D} B(\mu)(\zeta) =
\sup_{\zeta \in \partial \D} \int_{\D} P_z (\zeta) \, d \mu(z) =
\sup \left\{ \int_{\D} h  d\mu : h \in \Har_+(\D), h(0) =
1\right\},
\]
and for any positive harmonic function $h$,
\[
h(0) = \max_{\mu \in \mathcal B} \int_{\D} h  d\mu .
\]
\end{proposition}

\begin{proof}
Let $h=P[\nu]$, where $\nu\geq 0$ is a  measure on $\partial \D$. If $\mu$ has
balayage bounded by $c$,
\[
\int_{\D} h(z)  d\mu(z)=
     \int_{\partial \D} \int_{\D} P_z (\zeta) d\mu(z) d \nu (\zeta)
\le c  \nu (\partial \D) = c h(0).
\]
Conversely, since $z \mapsto P_z(\zeta)$ is a harmonic function for
any fixed $\zeta$, $\int_{\D} P_z (\zeta) \, d \mu(z)$ is
pointwise defined. Pick a sequence $\zeta_n$ such that
\[
\lim_{n\to\infty}\int_{\D} P_z (\zeta_n) \, d \mu(z) =
\sup_{\zeta \in \partial \D}  \int_{\D} P_z (\zeta) \, d \mu(z),
\]
where the supremum on the right hand side might a priori be infinite. Since the
set $E:=\{h\in \Har_+(\D):h(0)=1\}$ is uniformly bounded on compact sets in
$\D$, a normal family argument shows that $\sup\{\int_\D hd\mu:h\in
E\}<\infty$. Observe that the mapping $z\mapsto P_z(\zeta_n)$ is in $E$ for
every $\zeta_n$, $n\in \N$. Hence $\sup_n \int_\D P_z(\zeta_n)d\mu(z)<\infty$.

This proves that $\mu$ has bounded balayage, and the equalities between
constants that we had announced.
\end{proof}

The next result is a refined version of Theorem \ref {thmequivBB} stated
in the introduction.

\begin{theorem}\label{equivBB}
Let $\varphi$ be a nonnegative Borel  function on
the unit disk. Then there exists a (positive) harmonic function
$h$ such that $h (z) \ge \varphi(z)$ for any $z \in \D$ if and
only if
\begin{equation}
\label{constBB}
M_\varphi := \sup_{\mu \in \mathcal B} \int_{\D}
\varphi \, d\mu < \infty.
\end{equation}
Furthermore,
\[
M_\varphi = \inf \left\{ h(0) : h \in \Har (\D), h \ge \varphi
\right\}.
\]
\end{theorem}

That (\ref{constBB}) is necessary is clear from the above
considerations. In order to prove that it is sufficient, we will
reduce ourselves to a discrete version of it.
We will use the dyadic squares introduced in \eqref{dyadsquares}.
As in the previous
section, choose a point
$z_{n,k}$ in each cube, say
\begin{equation*}
z(Q_{n,k})=z_{n,k} := (1- 2^{-n})
\exp({2\pi}k 2^{-n}).
\end{equation*}

Observe that by
Harnack's inequality,
there exists a universal constant $K$ such
that : if $z, z'$ lie in the same Whitney square
$Q_{n,k}$ (as defined in \eqref{dyadsquares}),
then $K^{-1} P_{z'} (\zeta) \le P_{z} (\zeta) \le K
P_{z'} (\zeta)$, for any $\zeta \in \partial \D$.

\begin{lemma}
\label{discBB} The function $\varphi$ satisfies condition
(\ref{constBB})  if and only if there exists a constant
$M'_\varphi$ such that for any sequence of nonnegative
coefficients $\{c_{n,k}\}$ such that
\begin{equation}
\label{sumPK}
\sup_{\zeta \in \partial \D} \sum_{n,k} c_{n,k} P_{z_{n,k}}
(\zeta) \le 1 ,
\end{equation}
then
\begin{equation}
\label{sumsupphi}
     \sum_{n,k}  c_{n,k} \sup_{Q_{n,k}} \varphi  \le M'_\varphi .
\end{equation}
Furthermore, $ C^{-1} M_\varphi \le  M'_\varphi \le C M_\varphi$,
where $C>1$ is an absolute constant.
\end{lemma}

\begin{proof}[Proof of Lemma \ref{discBB}]

Pick $z_{n,k}^* \in Q_{n,k}$ such that
$\varphi(z_{n,k}^*) \ge ( \sup_{Q_{n,k}} \varphi)/2$
and define the measure $\mu := \sum_{n,k} c_{n,k} \delta_{z_{n,k}^*}$.
Then, if $\{c_{n,k}\}$ satisfies (\ref{sumPK}),
\[
\int_{\D} P_z (\zeta) \, d \mu(z) = \sum_{n,k} c_{n,k}
P_{z_{n,k}^*} (\zeta) \le K  \sum_{n,k} c_{n,k} P_{z_{n,k}}
(\zeta) \le K.
\]
So if $\varphi$ satisfies (\ref{constBB}),
\[
\sum_{n,k}  c_{n,k} \sup_{Q_{n,k}} \varphi  \le 2 \sum_{n,k}
c_{n,k} \varphi(z_{n,k}^*) = 2 \int_{\D} \varphi \, d\mu \le 2
K M_\varphi.
\]
The converse direction is easier, and left to the reader (it also
follows from the proof of the theorem, below).
\end{proof}

We now need a classical result in convex analysis.
Recall that the convex hull of a
subset $A\subset \mathbb R^d$ is defined as
\[
\Conv(A) := \bigl\{ \sum_{i=1}^N \a_i a_i : a_i \in
A, \a_i \ge 0, \sum_i \a_i = 1 \bigr\} .
\]
If we write $\mathbb R_+ A := \{ \la x : \la \ge 0, x \in A \}$,
then the conical convex hull of $A$ is defined as
\[
\Cone(A) := \Conv(\mathbb R_+ A) =
\bigl\{ \sum_{i=1}^N \a_i a_i : a_i \in A, \a_i \ge 0 \bigr\} .
\]
When $A$ is a finite set, the
conical convex hull is equal to its closure:
$\overline{\Cone} (A) = \Cone(A)$
(for this and other facts, see \cite{HULL}). The key fact for us
will be the generalized form of the Minkowski-Farkas Lemma
\cite[Chapter III, Theorem 4.3.4]{HULL} that we cite here only
for finite $A$. Let $\langle \cdot ,
\cdot \rangle$ stand for the standard Euclidean scalar product in
$\mathbb R^d $.

\begin{theorem}
\label{Farkas} Let $(a_j,b_j) \in \mathbb R^d \times \mathbb R$,
$1\le j \le N$, be such that
$X := \{  x \in \mathbb R^d : \langle a_j, x \rangle \le b_j
\} \neq \emptyset$.
Denote $A:= \{(a_j,b_j), 1\le j \le N\} \subset \mathbb R^d \times
\mathbb R$.
Then the following properties are equivalent for
$(v,r) \in  \mathbb R^d \times \mathbb R$:
\begin{itemize}
\item[(a)] For any  $x \in X$, $\langle v , x \rangle\le r$.
\item[(b)] $(v,r) \in \Cone (A)$.
\end{itemize}
\end{theorem}

We will use the following special case. For a vector  $v \in
\mathbb R^d$, the coordinates are denoted by $v^i$, $1\le i \le
d$. Also, $\mathbb R_+^d$ denotes the set of points of $\mathbb
R^d$ with nonnegative coordinates.

\begin{corollary}
\label{posFarkas} Given $a_j\in \mathbb R^d $, $1\le j \le N$,
let $X_+ := \{  x \in \mathbb R_+^d : \langle a_j, x \rangle
\le 1 \}$,
and suppose that $X_+ \neq \emptyset$. 
Then the following properties are equivalent for
$v \in  \mathbb R_+^d $:
\begin{itemize}
\item[(a)] For any $x \in X_+$, $\langle v , x \rangle \le 1$.
\item[(b)] There exist $\a_j \ge 0, 1 \le j \le N$
such that $\sum_{j=1}^N \a_j = 1$ and for any $i= 1, \dots, d$,
\[
v^i \le \sum_{j=1}^N \a_j a_j^i.
\]
\end{itemize}
\end{corollary}

\begin{proof}
Let $\{ e_i \}_{1\le i\le d}$ be the canonical basis of $\mathbb
R^d$ and consider
\[
A := \{ (a_j,1),  1\le j \le N \} \cup \{ (-e_i,0), 1\le i \le d
\} .
\]
Then $X_+$ corresponds to the $X$ in Theorem \ref{Farkas},
from what we see that (a) implies that there exist
$\a_j \ge 0,\be_i \ge 0$,  $1 \le j \le N$,  $1 \le i \le d$,
such that
\[
(v,1) = \sum_{j=1}^N \a_j (a_j,1) - \sum_{i=1}^d \be_i (e_i,0) .
\]
When applied to each coordinate, this yields $1=\sum_{j=1}^N \a_j$ and
\[
v^i = \sum_{j=1}^N \a_j a_j^i - \be_i \le \sum_{j=1}^N \a_j a_j^i.
\]
The converse implication is immediate.
\end{proof}

\begin{proof}[Proof of Theorem \ref{equivBB}]
Suppose that  $\varphi$ satisfies (\ref{constBB}). For each
nonnegative integer $m$,
we define
\[
a_j := \left( P_{z_{n,k}} (\exp( i j \cdot 2^{-m} 2\pi)) \right)_{
\begin{subarray}{l}
0\le n \le m\\
0 \le k \le 2^n -1
\end{subarray}} \mbox{ for } 0 \le j \le 2^m -1,
\]
$d:= \sum_{n=0}^m 2^n$ and
\begin{multline*}
X_+ := \bigl\{ \{ c_{n,k} \}_{
\begin{subarray}{l}
0 \le n \le m\\
0 \le k \le 2^n -1
\end{subarray}}
\in \mathbb R_+^d :   \\
\sum_{
\begin{subarray}{l}
0 \le n \le m\\
0 \le k \le 2^n -1
\end{subarray}}
c_{n,k} P_{z_{n,k}} (\exp(i j \cdot  2^{-m} 2\pi)) \le 1, \mbox{ for } 1
\le
j \le 2^m -1 \bigr\} .
\end{multline*}
Obviously, $X_+$ is not empty:  for instance $c_{0,0}=1$ and
$c_{n,k}=0$ for $n\ge 1$ gives a point in $X_+$. We claim that any
$\{ c_{n,k} \} \in X_+$ will satisfy (\ref{sumPK}) up to a
constant. Indeed, for any $\theta \in [0,2\pi)$, there is an index
$j<2^m$ so that $j \cdot 2^{-m} 2\pi \le \theta < (j+1) \cdot
2^{-m} 2\pi$, therefore by Harnack's inequality, for any $z$ such
that $|z|\le 1-2^{-m}$,
\[
P_z(\eit) = P_{z\exp(i(j\cdot 2^{-m}2\pi-\theta))} (\exp( i j
\cdot 2^{-m} 2\pi)) \le K P_z (\exp( i j \cdot 2^{-m} 2\pi).
\]
Therefore $\{ K^{-1} c_{n,k} \} $  satisfies (\ref{sumPK}),
and by Lemma \ref{discBB} and the hypothesis, $\varphi$ satisfies
(\ref{sumsupphi}) with constant $K M'_\varphi$. Corollary
\ref{posFarkas} then implies the existence of positive
coefficients $( \a_j^m )_{j=0}^{ 2^m-1}$ with sum
equal to $K M'_\varphi$,
such that
\[
\sup_{Q_{n,k}} \varphi \le \sum_{j=0}^{2^m-1} \a_j^m P_{z_{n,k}} (\exp(
i j \cdot 2^{-m} 2\pi)) = \int_{\partial\mathbb D} P_{z_{n,k}} d \nu^m,
\]
where $\nu^m$ is the discrete measure on the circle given by
the following combination of Dirac masses:
\[
\nu^m = \sum_{j=0}^{2^m-1} \a_j^m \delta_{\exp( i j \cdot 2^{-m}
2\pi)} .
\]
Since the mass of $\nu^m$ is uniformly bounded by $K
M'_\varphi$, we can take a weak* limit $\nu$ of this
sequence of measures, so that for any $(n,k)$,
\[
\sup_{Q_{n,k}} \varphi \le \int_{\partial\mathbb D} P_{z_{n,k}} d \nu =
h(z_{n,k}),
\]
where $h:=P[\nu]$. Harnack's
inequality now implies that there is an absolute constant $C_1$
such that $C_1 h(z) \ge \varphi (z)$ for any $z \in \D$. This
proves the theorem, with the inequality
\[
\inf \left\{ h(0) : h \in \Har (\D), h \ge \varphi \right\} \le C_1
K M'_\varphi \le C C_1 K M_\varphi .
\]
The constants $C$, $K$ and $C_1$ only depend on the
discretization we have chosen. Picking a
discretization with smaller ``squares", we may make all three
constants as close to $1$ as we wish.
\end{proof}

We finish this section with the proof of Corollary~\ref{decrease}.

\begin{proof}[Proof of Corollary~\ref{decrease}]
Given a non-Blaschke sequence $\Lambda$, arguing as in \cite{NPT} one
can show that
there exists a function $f\in H^\infty(\D)$ in the unit disc with
$|f(\lambda)|< \varepsilon_\lambda$ for any $\lambda \in \Lambda$ if and
only
if $\Lambda$ is the union of a Blaschke sequence and a sequence $\Gamma$
for
which there exists a positive harmonic function $h$ in the unit disc with
$h (\gamma) > - \log  \varepsilon_\gamma$ for all $  \gamma \in \Gamma$.
Then the result follows from Theorem~\ref{thmequivBB}.
\end{proof}

\section{Weaker conditions for the existence of
harmonic majorants}
\label{weakcond}

In this section we state first a sufficient condition implied by a result of
Borichev on a similar problem. On the other hand, we also  prove the necessary
condition of Proposition~\ref{strongneccond} and show that  it is not
sufficient.

The following was communicated to us by our colleague Alexander Borichev (see
\cite{NiTh} for a proof).

\begin{theorem}
\label{discrete}
Given a collection of nonnegative data $\{\varphi_{n,k}\} \subset \R_+$,
there
exists a finite positive measure $\nu$ on $\partial \D$ such that
\[
\frac{ \nu(I_{n,k})}{\sigma(I_{n,k})} \ge \varphi_{n,k}
\]
if and only if
\begin{equation}
\label{discNSC}
S=:\sup \bigl\{ \sum_{(n,k) \in A} \varphi_{n,k} \sigma(I_{n,k}) :
\{I_{n,k}\}_{ (n,k) \in A}  \mbox{ is a disjoint family}
\bigr\}
< \infty.
\end{equation}
\end{theorem}

This is an analogue of the discretized version of
Theorem~\ref{thmCNSN}(d), (as in Lemma~\ref{discBB}) obtained by
considering only measures of type $\mu_A := \sum_{(n,k) \in A}
\sigma(I_{n,k}) \delta_{z_{n,k}}$, and by replacing the Poisson kernel
$P_z$ by the ``square" kernels
\begin{equation*}
K_z (\eit) := K_{I_z} (\eit) :=\frac1{\sigma(I_z)} \chi_{I_z}(\eit).
\end{equation*}
Here $I_z$ denote the intervals defined in \eqref{shadow} and $\chi_E$
stands for the characteristic function of $E$.

The similarity of Theorem \ref{thmCNSN} with this result
leads us to an:

{\bf Open Question.}  Is condition (d) in Theorem \ref{thmCNSN} still
sufficient if we restrict it to $\{ c_\lambda \}$ such that for any
$\lambda \in \Lambda$, $c_\lambda =0$ or $(1-|\lambda|)$?

Theorem \ref{discrete} together with the estimate $K_z\lesssim P_z$
provide a sufficient (but not necessary) condition for domination
by true harmonic functions, which is clearly less restrictive than
requiring
that $M\varphi \in L^1(\partial \D)$, but easier to check in concrete
examples than the characterizing condition of Theorem
\ref{thmequivBB}.

\begin{corollary}
\label{compdiscrete}
Any positive function $\varphi$ such that $\varphi_{n,k}:= \sup_{Q_{n,k}}
\varphi$ satisfies \eqref{discNSC} admits a harmonic majorant.
On the other hand, the positive harmonic function $z \mapsto P_z (1)$ does
not satisfy  \eqref{discNSC} for certain choices of $A$.
\end{corollary}

\begin{proof}

It is well known and easy to see that 
there exists a constant $c$ such
that $P_z \ge c K_{I_{n,k}}$ for any $z \in Q_{n,k}$
(the constant $c$ depends on the aperture $\a$ of the Stolz angle).
Therefore, for any $z \in Q_{n,k}$
\[
P[\nu](z) \ge c \int_{\partial \D} K_{I_{n,k}} (\zeta) d\nu(\zeta)\geq
c\frac{\nu(I_{n,k})}{\sigma(I_{n,k})} \ge c\varphi_{n,k} = c\sup_{Q_{n,k}}
\varphi,
\]
which proves that $P[(1/c)\nu]$ is the harmonic majorant we are looking
for.

To see that the condition
is not necessary,  consider any  $A \subset \{ (n,1) : n \in \mathbb
N\}$. Then the intervals $I_{n,1}$ are all disjoint; however
$P_{z_{n,1}}(1) \simeq  2^n \simeq \sigma(I_{n,1})^{-1}$, so that
condition \eqref{discNSC} will fail (the sum is comparable to $\# A$).
\end{proof}

In the same way as in Corollary \ref{CNgeom}, Corollary
\ref{compdiscrete} and Proposition \ref{strongneccond} imply the
following result.  For $Q=Q_{n,k}$, write $I(Q)=I_{n,k}$ (the radial
projection of the square to an arc of the circle).

\begin{corollary}
\label{coroBori}
Assume that $\Lambda$ is contained in a
union $A$ of Whitney squares $Q$ of center $z(Q)$ and that
\[
\sup \bigl\{
\sum_{Q \in A'} (1 - |z(Q)|) \sup_{\lambda \in \Lambda \cap Q}
\log|B_{\lambda} ( \lambda)|^{-1}
\bigr\} < \infty,
\]
where the supremum is taken over all
$A' \subset A$ such that $\{I(Q), Q \in A'\}$
is a disjoint family,
then $\Lambda$ is interpolating for the Nevanlinna class.
\end{corollary}

We move next to the proof of the necessary condition in terms of the
Hardy-Littlewood maximal function.

\begin{proof}[Proof of Proposition \ref{strongneccond}]
(a) The problem can be localized, so we may work on the upper half plane,
$\mathbb{C}_+ := \{ x + i y : y > 0 \}$, with $I_{x+iy} := (x-y,x+y)$,
restricting ourselves to positive harmonic functions which are
Poisson integrals of positive measures with finite mass. Here the Poisson
kernel is given by
\[
P_{x+iy}(s) = \frac1\pi \frac{y}{(x-s)^2+y^2}.
\]
For a measurable
set $E \subset \mathbb R$ write $|E|=\sigma(E)$.
Also, we only need to look at boundary points in a fixed bounded interval,
say $-1 \le x \le 1$.

For any $t>0$,
let $E_t := \{ s \in [-1,1] : \varphi^H (s) > t \}$. For any $s \in
E_t$,
there exists $z=z(s)$ and $J=J(s) \supset I_z$ such that
\begin{equation}
\label{maxint}
\varphi(z) \int_J \chi_{I_z} > t |J|,\quad \mbox{ i.e. }\quad
\varphi(z) |I_z| > t |J|.
\end{equation}
By Vitali's
covering lemma, there exist an absolute constant $c_1 \in (0,1)$ and a
disjoint family of intervals $J_j := J(s_j)$, $1\le j \le N$, such that
$\sum_j  |J_j | \ge c_1 |E_t|$.

Write $z_j:= z(s_j)=: x_j+iy_j$. Note that since the point $z_j$ is
contained in the ``tent" over $I_{z_j}$ (therefore in the tent over $J_j$)
the points $z_j$ are separated in the hyperbolic metric.

Now define new points $z'_j$ in the following way : let
$y'_j :=   |J_j|/2 = y_j |J_j|/|I_{z_j}| \ge y_j$ and
$z'_j := x_j + i y'_j$.
Note that $|J_j \cap I_{z'_j}| \ge  |J_j|/2$.

We claim that $h(z'_j) \ge  t$, where $h$ is a harmonic majorant
of $\varphi$. Indeed, writing $h = P[\mu]$,
$$
h(z'_j) = \frac1{ \pi y'_j}
\int_{\R} \frac1{1+ \bigl( \frac{t-x_j}{y'_j} \bigr)^2} \, d\mu(t)
\ge
\frac1{\pi y'_j}
\int_{\R} \frac1{1+ \bigl( \frac{t-x_j}{y_j} \bigr)^2} \, d\mu(t)
= \frac{y_j}{y'_j} h(z_j),
$$
and, by (\ref{maxint}),
    $h(z_j) \ge \varphi(z_j) > t |J_j|/|I_{z_j}| =  t y'_j/y_j$.

Therefore, since $Mh \in L^1_w(\mathbb R)$,
\[
\frac{c_1}2 |E_t| \le \frac12 \sum_j  |J_j | \le
\sum_j |J_j \cap I_{z'_j}| \le  \left| \{ M h > t \} \right| \le
\frac{C_h}t.
\]

(b) Similarly.
\end{proof}

We now give two examples showing that the necessary condition of
Proposition \ref{strongneccond} is strictly stronger than that
of Proposition \ref{maxfcnharm} but still not sufficient.

\begin{proposition}
\label{notsuffstronger}
\begin{itemize}
\item[(a)]
There are functions $\varphi$ such that $\varphi^H \in L^1_w(\partial \D)$,
but that do not admit a harmonic majorant.
\item[(b)]
There are functions $\varphi$ such that $M \varphi \in L^1_w(\partial \D)$,
but $\varphi^H \notin L^1_w(\partial \D)$.
\end{itemize}
\end{proposition}

\begin{proof}
The proof  rests on the following family of examples.  Note that it is easy
to turn those examples into examples of sequences which are (or are not)
interpolating for the Nevanlinna class.

Again we will work on $\mathbb{C}_+$.
Our functions $\varphi$ will vanish everywhere on the upper half plane,
except on the sequence $\lambda_k := x_k + i y_k$, where
$x_k=k^{-\alpha}$ and $y_k=k^{-\beta}$. To
ensure that $y_k \le  (x_{k+1}-x_k)^2$ we take $\beta\geq 2(\alpha+1)$.
With this choice, it can be deduced from Proposition~\ref{thmnafta2}
(or the remark before
Lemma \ref{superseparat}),
that a  necessary and sufficient condition for
the existence of a harmonic majorant is  that $M \varphi \in L^1$,
that is,
\begin{equation}\label{CNSvarphi}
   \sum_k \varphi(\lambda_k) y_k < \infty.
\end{equation}

We note that
\[
\chi_{I_{\lambda_k}}^* (x) = \frac2{1+\max(1,\frac{|x-x_k|}{y_k})}.
\]

Henceforth we only study data $\{\varphi_k\} := \{\varphi(\lambda_k)\}$
which are increasing sequences of
positive numbers tending to infinity. We also assume that
$\{(\varphi_k y_k + \varphi_{k+1} y_{k+1})/(x_{k}-x_{k+1}) \}_k$
forms an increasing sequence. Let
$k_0 (t) := \min \{ k : t < \varphi_k \}$.
The necessary condition arising from the fact that
$M \varphi \in L^1_w (\mathbb R)$ reads
\begin{equation}\label{maxfex}
\sum_{k \ge k_0(t)} y_k \simeq k_0^{1-\beta}(t)\le \frac{C}t, \qquad
\forall t >0.
\end{equation}

This condition will be assumed for both examples.

Now, for $k \ge k_0 (t)$, define
$J_k := \{ x : \varphi_k \chi_k^* (x) > t \}
\simeq ( x_k- y_k \varphi_k/t , x_k + y_k \varphi_k/t )$,
and let $k_1 (t) := \min \{ k : J_k \cap J_{k+1} \neq \emptyset \}$.
Then,
\[
\bigcup_{k \ge k_1(t)} J_k = (0, x_{k_1(t)} + y_{k_1(t)}
\frac{\varphi_{k_1(t)}}{t} )=(0,k_1^{-\alpha}(t)+ k_1^{-\beta}(t)
\frac{\varphi_{k_1(t)}}{t})
\]
and
\begin{equation}
\label{measgtl}
\left| \{ x : \varphi^H (x) > t \} \right| \simeq
   k_1^{-\alpha}(t) + k_1^{-\beta}(t)
   \frac{\varphi_{k_1(t)}}t + \frac2t \sum_{k=k_0(t)}^{k_1(t)-1}
   \frac{\varphi_k}{k^{\beta}} \simeq
k_1^{-\alpha}(t)  + \frac2t \sum_{k=k_0(t)}^{k_1(t)}
\frac{\varphi_k}{k^{\beta}} .
\end{equation}

In order to prove (a), choose $\varphi_k := \eps_k k^{\beta-1}$. Since
$t\simeq \eps_{k_0(t)} k_0(t)^{\beta-1}$,
condition \eqref{maxfex} becomes
that $(\eps_k)_k$ remains bounded above, while the necessary and
sufficient condition (see \eqref{CNSvarphi}) is
$$
\sum_k \frac{\eps_k}k < \infty.
$$

With $\eps_k:= (\log k)^{-1}$, this
condition fails, so that $\varphi$ admits no harmonic
majorant.

However,
$k_0 (t) \simeq (t \log t )^{1/(\beta-1)}$.
Since $x_{k}-x_{k+1} \simeq k^{-\alpha-1}$, then
$1/k_1(t)^{\alpha+1} \simeq \eps_{k_1(t)}/(t k)$,
thus
$k_1(t)  \simeq ( t/\eps_{k_1(t)})^{1/\alpha}$,
and $k_1 (t) \simeq (t \log t )^{1/\alpha}$.

Therefore equation (\ref{measgtl}) becomes
\begin{align*}
\left| \{ x : \varphi^H (x) > t \} \right|& \simeq
\frac1{t \log t} + \frac2t \sum_{k=k_0(t)}^{k_1(t)} \frac1{k \log k}
\simeq \frac1{t \log t} +
\frac2t \log \left( \frac{\log k_1(t)}{\log k_0(t)} \right)\\
&\le \frac1{t \log t} + \frac{C}t \le \frac{C'}t,
\end{align*}
and this choice of $\varphi$ does satisfy the necessary condition given
in Proposition \ref{strongneccond}.

To prove the second statement in
the Lemma, choose $\eps_k:= 1$. With similar
but easier calculations one sees that
$k_0 (t) \simeq t^{1/(\beta-1)}$ and
$k_1 (t) \simeq t^{1/\alpha}$.
Therefore (\ref{measgtl}) becomes
$$
\left| \{ x : \varphi^H (x) > t \} \right| \simeq
\frac1{t } + \sum_{k=k_0(t)}^{k_1(t)-1} \frac1{k }
\simeq
\frac1{t } +\frac2t \log \left( \frac{ k_1(t)}{ k_0(t)} \right)
\simeq \frac{\log t}t,
$$
so the weak $L^1$ condition fails for $\varphi^H$, even though $\varphi$
satisfies the necessary condition in Proposition \ref{maxfcnharm}.
\end{proof}

\section{Proofs of the geometric conditions}
\label{GeomCond}

\begin{proof}[Proof of Corollary \ref{CNgeom}]

Since
\[
I_\lambda\subset \bigl\{\zeta\in\T : M\varphi_\Lambda(\zeta)\geq \log
|B_\lambda(\lambda)|^{-1}\bigr\} ,\quad \lambda\in\Lambda,
\]
to prove (a) and (b) it suffices to apply condition (a) of Corollary
\ref{maxfct}. Statement (c) follows from the next  Lemma applied to
$\varphi_\Lambda$.

\begin{lemma}
\label{L1suffcond}
Let $\varphi : \Dd \longrightarrow \mathbb R_+$ satisfy
$\sum_{\lambda \in \Dd} (1-|\la|) \varphi(\la) < \infty$. Then $\varphi$
admits a harmonic majorant.
\end{lemma}

\begin{proof}
Let $u=\sum_\lambda
\varphi(\lambda) \chi_{I_{\lambda}}$.  By assumption  $u\in
L^1(\T)$ and obviously $M\varphi \leq u$, hence the result
follows from Corollary \ref{maxfct}(b).
\end{proof}

Parts (b) and (c) also follow directly from Theorem~\ref{thmCNSN}(d), by a
simple argument based on the $\ell^1$, $\ell^\infty$ duality.
\end{proof}

\begin{proof}[Proof of Proposition \ref{thmnafta}]

It is enough to consider the case where $\Lambda$ is contained in
only one Stolz angle. Indeed, if
$\Lambda=\bigcup_{i=1}^n\Lambda_i$ with $\Lambda_l\subset
\Gamma_{\zeta_l}$, $l=1,\ldots,n$, and $\zeta_i\neq\zeta_j$, then
\[
\lim\limits_{z\to \zeta_i, z\in \Gamma_{\zeta_i}}
|B_{\Lambda_j}(z)|=1, \quad j\neq i, 
\]
so that $\log
|B_{\lambda}(\lambda)|^{-1}$ behaves asymptotically like $\log
|B_{\Lambda_i\setminus\{\lambda\}}(\lambda)|^{-1}$ in $\Gamma_{\zeta_i}$
(here $\lambda\in \Lambda_i$). Also, we can assume that the
sequence is radial (this means that we replace the initial
sequence by one which is in a uniform pseudo-hyperbolic
neighborhood of the initial one; by Harnack's inequality such a
perturbation does not change substantially the behavior of
positive harmonic functions).

According to Corollary \ref{CNgeom} it is enough to prove the sufficiency of
the conditions. Let us first show that condition \eqref{CN} implies
interpolation in $N^+$. In order to construct a function $w\in L^1(\T)$ meeting
the requirement of Theorem \ref{thmCNS}(c) assume that
$\Lambda=\{\lambda_n\}_n\subset [0,1)$ is arranged in increasing order and set
$\tilde{\eps}_n=(1-|\lambda_n|)\log |B_{\lambda_n} (\lambda_n)|^{-1}$. Clearly
there exists a decreasing sequence $(\eps_n)_n$ with $\tilde{\eps}_n\le\eps_n$,
$n\in \N$, and $\lim_n \eps_n=0$. Now, if $I_n=I_{\lambda_n}$,  set
$J_n=I_n\setminus I_{n+1}$, $\beta_n= \eps_n-\eps_{n+1}$, and set
\beqa
    w(\zeta)=\sum_n \frac{\beta_n}{\sigma(J_n)} \chi_{J_n}(\zeta),\quad
    \zeta\in\T.
\eeqa
Then $w\in L^1(\T)$, and
\begin{eqnarray*}
    P[w](\lambda_n) &\ge&  \int_{I_n} P(\lambda_n,\zeta)
      \sum_k \frac{\beta_k}{\sigma(J_n)}
     \chi_{J_k}(\zeta) \,d\sigma(\zeta)
    \gtrsim \sum_{k\ge n}\frac{\beta_k}{\sigma(J_n)}
      \frac{1}{(1-|\lambda_n|)}\int_{J_k}\,d\sigma(\zeta)\\
    &=&\frac{\sum_{k\ge n} \beta_k}{1-|\lambda_n|}
    =\frac{\eps_n}{1-|\lambda_n|}
     \ge \frac{\tilde{\eps}_n}{1-|\lambda_n|}
=\log   |B_{\lambda_n} (\lambda_n)|^{-1}.
\end{eqnarray*}
This and Theorem \ref{thmCNS} prove the assertion.

The proof for the Nevanlinna class is even simpler.
Set $d\mu_s=  \delta_{1}$, the Dirac mass on $1\in\T$. From  \eqref{CNN}
we get
\[
    \log  |B_{\lambda_n}(\lambda_n)|^{-1}\lesssim \frac{1}{1-|\lambda_n|}
    \lesssim  P[\mu_s](\lambda_n),
\]
and we finish by applying Theorem \ref{thmCNSN}.
\end{proof}

\begin{proof*}{{\it Proof of Proposition \ref{supersep}.}}
By Corollary \ref{CNgeom}(c), we already know that \eqref{CS} is a
sufficient condition for $\Lambda$ to be interpolating for $N^+$.
Conversely, suppose that $\Lambda$ is interpolating for $N$, that is, 
$\varphi_\Lambda$ admits a harmonic majorant. Since $\mu_\Lambda$ has
bounded balayage,  then  
$\int_{\Dd} \varphi_\Lambda d \mu_\Lambda
< \infty$, which is exactly \eqref{CS}.
\end{proof*}

\begin{proof*}{{\it Proof of Proposition \ref{embedyana}.}}
It is obvious that (a) implies (b). If we assume (c), $\mu_\Lambda$ will act
against any positive harmonic function. Suppose $F \in N$. As seen in 
Section \ref{trace}, there exists a positive harmonic function $h$
so that $\log^+|F|\le h$. Thus, taking $\mu_{\lambda}$ as in 
\eqref{eqnp7},
\[
\sum_{\la \in \Lambda} (1-|\la|) \log^+|F(\la)|
=
\int_{\Dd} \log^+|F(\la)| \, d\mu_\Lambda (\la) \le \int_{\Dd} h(\la) \,
d\mu_\Lambda (\la) < \infty.
\]
Finally, to prove that (b) implies (c), suppose that (c) doesn't hold, i.e.
$g_\Lambda := \sum_\la (1-|\la|) P_\la$ is unbounded. Since $g_\Lambda$ is
lower semi-continuous, this implies that $g_\Lambda\notin
L^{\infty}(\partial\D)$. Since $L^\infty$ is the dual of $L^1$, 
there exists $f\in L^1(\T)$ such that $\int_{\T} fg_\Lambda
=\infty$. Taking an outer function $F \in N^+$ with $\log|F|=P[f]$ we see
that 
\[
\sum_{\la \in \Lambda} (1-|\la|) \log|F(\la)| = \sum_{\la
\in \Lambda} (1-|\la|) \int_{\T} P_\lambda f 
= \int_{\T} f g_\Lambda
=\infty,
\]
so (b) doesn't hold. 
\end{proof*}

\begin{proof*}{{\it Proof of Naftalevi\v c's theorem.}} Assume that $\Lambda$
is contained in a finite union of Stolz angles and \eqref{CNN} holds. By
Proposition \ref{thmnafta}, $\Lambda\in\Int N$, hence the trace $N|\Lambda$ is
given by the majorization condition of Theorem~\ref{thmCNSN}(b). Taking as
majorizing function the Poisson integral of the sum of the Dirac masses at the
vertices, we see that $N|\Lambda\supset l_{\text{Na}}$.

Conversely, if $N\vert \Lambda = l_{\text{Na}}$ then the
trace is  ideal, so $\Lambda$ is free interpolating and by Corollary
\ref{CNgeom}(b) \eqref{CNN} holds.
According to Theorem~\ref{thmCNSN}(b) and the definition of $l_{\text{Na}}$,
the function 
\[
\varphi(z)= 
\begin{cases} (1-|\la|)^{-1}\quad&\textrm{if $z=\la \in \Lambda$}\\
\ 0\quad&\textrm{if $z\notin\Lambda$}
\end{cases}
\]
admits a harmonic majorant $h$. Let $h(z)=P[\mu](z)$ and consider the intervals
\[
I_z^\alpha=\{\zeta\in\partial\D : z\in \Gamma_\alpha(\zeta)\} .
\]
There exist constants $\alpha $ and $C_\alpha$ such that $\mu( I_{ z}^\alpha) >
C_\alpha$ for any  $z$ such that $h(z) \ge (1-|z|)^{-1}$.

If $\Lambda$ is not contained in a finite union of Stolz angles, then there is
an accumulation point $\zeta\in\T$ of $\Lambda' \subset \Lambda$ such that
$\Lambda'\not\subset \Gamma_\beta(\zeta)$ for any $\beta$. Pick $\beta >
\alpha$; then for $\lambda' \in \Lambda'$, $I_{\lambda'}^\alpha \not\ni \zeta$
and we can construct an infinite subsequence $\Lambda'' \subset \Lambda'$ such
that the Privalov shadows $\{ I_{\lambda'}^\alpha\}_{ \lambda' \in \Lambda''}$
are disjoint.  This prevents $h$ from being the Poisson integral of a finite
positive measure. \end{proof*}

We now give an example of a concrete separation condition implying that
$\mu_\Lambda$ has bounded balayage.

\begin{proposition}
\label{thmnafta2}
Assume that $\Lambda\subset\D$ is contained in the union of a
family $A$ of Whitney squares such that
\[
|\Arg (z(Q)) - \Arg(z(L))| > g^{-1}(1-|z(Q)|)
\]
for any $ Q , L \in A$, $Q\neq L$,
where  $z(Q)$ is the center of $Q$ and
$g$ is a positive function,  with $g(x)/x$ decreasing and 
\[
\int_0 \frac{g(x)}{x^2} < \infty \, .
\]
Then $\Lambda\in Int N$ if and only if $\Lambda\in Int N^+$, 
and if and only if
\[
\sum_{Q \in A} (1 - |z(Q)|) \sup_{\lambda \in \Lambda \cap Q}
\log|B_{\lambda} ( \lambda)|^{-1} < \infty
\]
\end{proposition}

Note that this covers some cases where $\mu_\Lambda$ has not bounded balayage,
even though another measure associated with the sequence will be (see the
proof).

In order to prove Proposition
\ref{thmnafta2} consider the ``Carleson window" $Q(\eit, r)$  centered at
$\eit$, of side $r$:
$$
Q(\eit, r):= \{ z \in \D : 1-|z|\le r, |\Arg (z) - \theta| \le
r \}.
$$
The next result is a Carleson-type condition which implies boundedness
of the balayage.

\begin{lemma}
\label{sufcond} Suppose that $\mu (Q(\eit,r))\le g(r)$, where $g$
is a nondecreasing function on $[0,2)$ with
$$
\int_0 \frac{g(x)}{x^2} dx < \infty.
$$
Then $\mu$ is a measure with bounded balayage.
\end{lemma}

A discrete version of this condition is
$$
\sum_n 2^{n} \sup_{\theta \in \mathbb R} \mu(Q(\eit,2^{-n})) <
\infty,
$$
as can be checked by writing a Riemann sum.

\begin{proof}
For any $t>0$, let $\Omega_t (\theta) := \{ z \in \D :  P_z (\eit)
\ge t\}$. This is a disk, tangent to the unit circle at the point
$\eit$, of radius $1/(t+1)$. Therefore $\Omega_t (\theta) \subset
Q(\eit,C/t)$ for $t\ge 1$, say, with $C>0$ an absolute constant.

Using the distribution function $\mu(\Omega_t(\theta))$ and the
fact that the measure $\mu$ is bounded,
we get the following estimate for the
balay\'ee of $\mu$:
\begin{eqnarray*}
\int_{\D} P_z (\eit) d\mu(z)&=& \int_0^\infty
\mu(\Omega_t (\theta)) dt\le
C_1 + \int_1^\infty \mu(\Omega_t (\theta)) dt
\le C_1 + \int_1^\infty \mu(Q(\eit,C/t)) dt \\
   &\le& C_1 + \int_1^\infty g(C/t) dt \le C_1 + C \int_0^1
\frac{g(x)}{x^2} dx <
\infty.
\end{eqnarray*}
\end{proof}

We will now compare measures satisfying the condition in Lemma
\ref{sufcond}, measures with bounded balayage and Carleson measures.
Each set is included in the next, and the examples will show that
both inclusions are strict.

{\bf Example 1.}  Let $\a = \{\a_n\}$ be a sequence of nonnegative
reals.  Let $\mu_\a$ be the measure concentrated on the circles
centered at the origin of radius $1-2^{-n}$ given in dual terms by
$$
\int_{\mathbb D} f (z) d \mu_\a (z) := \sum_{n\geq 1} \a_n \frac1{2\pi}
\int_0^{2\pi} f((1-2^{-n})\eit) d \theta .
$$

One can check that $\mu_\a$ is a Carleson measure if and only if it
has bounded balayage and this happens if and only if $\sum_n \a_n
< \infty$. Also $\mu_\a$ satisfies the condition in Lemma
\ref{sufcond}  if and only if $\sum_n \sum_{k\ge n} \a_k <
\infty$.

{\bf Example 2.} Let $m$ be a nonnegative-valued function on the
interval $[0,1)$.  Let $\mu_m$ be the measure concentrated on the
ray from the origin to $1$ given by
$$
\int_{\mathbb D} f (z) d \mu_m (z) := \int_0^1 f(x) m(x) dx .
$$

One can check that $\mu_m$ is a Carleson measure if and only if
there exists a constant $K$ such that
$$
\int_{1-\delta}^1 m(x) dx \le K \delta, \qquad \forall\delta >0
$$
and  $\mu_m$ is a measure with bounded balayage if and only if it
satisfies the condition in Lemma \ref{sufcond}, which happens if
and only if
$$
\int^1 \frac{m(x)}{1-x} dx < \infty.
$$

In particular, if we take $\a_k = k^{-\gamma}$, with $1<\gamma \le
2$, $\mu_m$ is a measure with bounded balayage but it does not satisfy
the condition in Lemma \ref{sufcond}; if we take $m(x) = 1$,
$\mu_m$ is a Carleson measure, but it does not have bounded balayage.

In view of Proposition \ref{supersep}, among other things, it is interesting
to understand for which separated sequences $\Lambda$ the corresponding
measure $\mu_\Lambda$ has bounded balayage. It is easy to see that this is the
case when $| \lambda/|\lambda| - \lambda'/|\lambda'| | \ge
(1-|\lambda|)^{1/2}$, $\lambda'\neq \lambda$, but more is true.

\begin{lemma}
\label{superseparat} Suppose that $g$ is a positive valued
function such that $g(x)/x$ is increasing and
$$
\int_0 \frac{g(x)}{x^2} dx < \infty.
$$
Let $g^{-1}$ stand for the inverse function of $g$. Then, if
we have a sequence $\Lambda\subset\D$  such that
$$
\left| \frac{\lambda}{|\lambda|} - \frac{\lambda'}{|\lambda'|} \right| \ge
g^{-1}(1-|\lambda|), \quad \forall \lambda'\neq \lambda,
$$
the measure $\mu_\Lambda$ has bounded balayage.
\end{lemma}

Examples of such functions $g$ are given by $x
(\log\frac1x)^{-1-\eps}$, with $\eps >0$. In that case, $g^{-1}(t)
\simeq t (\log\frac1t)^{1+\eps}$.

On the other hand, we can see that for the above lemma to hold,
we must have $g^{-1}(t) >> t$. More precisely,
take the sequence in the upper half-plane given by
$$
\lambda_k := e^{-k} + i k^{-1/2} e^{-k} .
$$
Then, $\Re\; \lambda_k- \Re\; \lambda_{k+1} \simeq  e^{-k}$,
so the sequence $\{\lambda_k\}_k$ verifies the separation condition in
Lemma
\ref{superseparat} with $g^{-1}(x) \simeq x (\log \frac1x)^{1/2}$, but
$$
\sum_k (\Im\; \lambda_k) P_{\lambda_k} (0) \simeq
\sum_k \left( \frac{\Re\; \lambda_k}{\Im\;
\lambda_k}\right)^2 = \sum_k \frac1k = \infty.
$$

\begin{proof}[Proof of Lemma \ref{superseparat}]
Let $\theta \in [0,2\pi)$. By hypothesis, there is
at most one $\liL$ such that
$$
\theta \in J_\lambda := \bigl( \Arg (\lambda) - \frac12
g^{-1}(1-|\lambda|) , \Arg (\lambda) + \frac12 g^{-1}(1-|\lambda|)
\bigr).
$$
Let $\mu':= \sum_{\lambda'\neq \lambda} (1-|\lambda'|) \delta_{\lambda'}$.
Then
$$
\int_{\mathbb D} P_z (\eit) d \mu_\Lambda (z) =(1-|\lambda|) P_{\lambda}
(\eit)
+ \int_{\mathbb D} P_z (\eit) d \mu' (z) \leq C + \int_{\mathbb D}
P_z (\eit) d \mu' (z).
$$
By the proof of Lemma \ref{sufcond} for this specific
value of $\theta$, we see that it will be enough to check that for
some absolute constants $C_1 , C_2$, one has
$$
\mu' (Q(\eit, r)) \le C_1 g(C_2 r), \mbox{ for } 0<r<2.
$$
Consider $\Sigma_r := \{ \lambda' \neq \lambda : \lambda' \in Q(\eit, r)
\}$.
For any $\lambda' \in \Sigma_r$ we have
$$
\sigma(J_{\lambda'}) = g^{-1}(1-|\lambda'|) \le 2 |\theta - \Arg
(\lambda')|
\le 2 r,
$$
so the intervals $J_{\lambda'}$ are all contained in $[\theta - 3r, \theta
+ 3r]$. Since they are disjoint, $\sum\limits_{\lambda'\in \Sigma_r}
\sigma(J_{\lambda'}) \le 6r$.

Using that $g(x)/x$ is increasing we have
$$
\sup_{\lambda'\in \Sigma_r} \frac{1-|\lambda'|}{\sigma(J_{\lambda'})} \le
\frac{g\bigl(\sup\limits_{\lambda'\in \Sigma_r}
\sigma(J_{\lambda'})\bigr)}{\sup\limits_{\lambda'\in \Sigma_r}
\sigma(J_{\lambda'})}
\le \frac{g(2r)}{2r} .
$$
Finally,
$$
\mu' (Q(\eit, r)) = \sum_{\lambda'\in \Sigma_r} 1-|\lambda'| \le
\sup_{\lambda'\in \Sigma_r}
\frac{1-|\lambda'|}{\sigma(J_{\lambda'})} \sum_{\lambda'\in \Sigma_r}
\sigma(J_{\lambda'}) \le \frac{g(2r)}{2r}
6r = 4 g(2r).
$$
\end{proof}

\begin{proof}[Proof of Proposition \ref{thmnafta2}]

For each Whitney square $Q$ in $A$, let $\lambda (Q)$ be the point
in $\Lambda \cap Q$ such that
\[
\log|B_{\lambda(Q)} (\lambda(Q))|^{-1} = \max \{\log|B_{\lambda}
(\lambda)|^{-1} : \lambda \in \Lambda \cap Q \} .
\]
Let $\Sigma$ be the
sequence formed by $\{\lambda(Q) : Q \in A \}$. By Lemma
\ref{superseparat} the corresponding measure $\mu_\Sigma$ has bounded
balayage. Therefore, there exists a positive harmonic function $h$
with $h(\lambda(Q)) \ge \log|B_{\lambda(Q)} (\lambda(Q))|^{-1}$ if
and only if
$$
\sum ( 1- |\lambda(Q)|) \log|B_{\lambda(Q)} (\lambda(Q))|^{-1}  <
\infty  \, .
$$
According to condition (c) in Theorem \ref{thmCNSN} one deduces that
$\Lambda\in\Int N$  if and only if the
last sum converges. Furthermore, when this is the case, the function $h$
can always be taken quasi-bounded (see Lemma \ref{L1suffcond}), so that
interpolation can actually be performed in the Smirnov class.
\end{proof}

\section{Hardy-Orlicz classes}\label{nevanlinna}

Let $\phi:{\R}\lra [0,\infty)$ be a convex, nondecreasing
function satisfying
\begin{itemize}
\item[(i)] $\lim_{t\to \infty} \phi (t)/t =\infty$
\item[(ii)] $\Delta_2$-condition: $\phi (t+2) \le M \phi (t) +K$,
$t \ge t_0$ for some constants $M,K \ge 0$ and $t_0\in \R$.
\end{itemize}
Such a function is called strongly convex (see \cite{RR}),
and one can associate with it the corresponding \emph{Hardy-Orlicz
class}
\beqa
    \hf =\{f\in N^+:\int_\T\phi (\log |f(\zeta)|)\,d\sigma(\zeta)
    <\infty\},
\eeqa
where $f(\zeta)$ is the non-tangential boundary value of $f$ at
$\zeta\in {\T}$, which exists almost everywhere.
In \cite{Har}, the following result was proved.

\begin{theorem*}
Let $\phi$ be a strongly convex function satisfying (i), (ii) and
the $V_2$-condition:
\beqa
    2\phi(t)\le \phi(t+\alpha), \quad t\ge t_1
\eeqa
where $\alpha>0$ is a suitable constant and $t_1\in\R$. Then $\Lambda
\subset {\D}$ is free interpolating for $\hf$
if and only if
$\Lambda$ is a Carleson sequence, and in this case
\beqa
    \hf\vert \Lambda=\{a=(a_\lambda)_\lambda: |a|_{\varphi}=
    \sum_{\liL} (1-|\lambda|)
    \varphi(\log |a_\lambda|)<\infty\}.
\eeqa
\end{theorem*}

The conditions on $\phi$ imply that for all
$\hf$ there exist $p,q\in (0,\infty)$ such that $H^p\subset \hf\subset
H^q$.
In particular, the $V_2$-condition
implies the inclusion $H^p\subset \hf$ for some $p>0$. This
$V_2$-condition has
a strong topological impact on the spaces. In fact, it guarantees that
metric
bounded sets are also bounded in the topology of the space (and so the
usual
functional analysis tools still apply in this situation; see
\cite{Har} for more on this and for further references).
It was not clear whether this was only a technical problem or if there
existed a critical growth for $\phi$ (below exponential growth
$\phi(t)=e^{pt}$ corresponding to $H^p$ spaces) giving a breakpoint in
the behavior of interpolating sequences for $\hf$.

We can now affirm that this behavior in fact changes between exponential
and polynomial growth. Let $\phi$ be a strongly convex function
with associated Hardy-Orlicz space $\hf$. Assume moreover that
$\phi$ satisfies
\bea\label{deltaineq}
    \phi(a+b)\le c(\phi(a)+\phi(b)),
\eea
for some fixed constant $c\ge 1$ and for all $a,b\ge t_0$.
The standard example in this setting is $\phi_p(t)=t^p$ for $p>1$.
We have the following result.

\begin{theorem}
Let $\phi:\R\lra [0,\infty)$ be a strongly convex function such that
\eqref{deltaineq} holds. If there exists a positive weight $w\in
L^1(\T)$ such that
$\phi \circ w\in L^1(\T)$ and $\varphi_\Lambda\le P[w]$,
then $\Lambda\in \Int \hf$.
\end{theorem}

\begin{proof}
Note first that \eqref{deltaineq} implies that $\hf$
is an algebra contained in $N^+$, hence it is sufficient to interpolate
bounded sequences (see Remark~\ref{l-infinit}).
As in Section~\ref{condsuff}, we set
\beqa
    g(z)=\int_\T \frac{\zeta+z}{\zeta-z}
     w(\zeta)d\sigma(\zeta).
\eeqa
The reasoning carried out in
Section \ref{condsuff}
leads to an interpolating function of the form $fH\exp(g)$, with
$f\in \Hi$, and $H=(2+g)^2$
outer in $H^p$ for all $p<1$ (note that the measure $\mu$ defining $g$ here
is absolutely continuous,
in fact $\mu=w\,d\sigma$). Also, $H^p\subset \hf$ for any $p>0$ by
our conditions
on $\phi$. By construction,
$\int \phi (\log |\exp g|)=\int \phi \circ w <\infty$ so that $\exp(g)\in
\hf$.
Since $\hf$ is an algebra, we deduce that $fH\exp(g)\in\hf$.
\end{proof}

\begin{example}
We give an example of an interpolating sequence for
$\hf$ which is not Carleson,
thus justifying our claim that there
is a breakpoint between Hardy-Orlicz spaces verifying the $V_2$-condition
and those that do not.

Consider the functions $\phi_p$ and let
$\Lambda_0=\{\lambda_n\}_n\subset \D$ be
a Carleson sequence  verifying $I_n\cap I_k=\emptyset$, $n\neq k$, where
$I_n$ are the arcs defined in \eqref{shadow}. Since
$\sum_n (1-|\lambda_n|)<\infty$, there exists a strictly increasing
sequence of
positive numbers $(\gamma_n)_n$ such that  $\sum_n
(1-|\lambda_n|)\gamma_n<\infty$ and $\lim_{n\to\infty}\gamma_n=\infty$.
Setting 
\beqa
    w=\sum_n \gamma_n^{1/p} \chi_{I_n},
\eeqa
we obtain $\int \phi_p \circ w
=\sum_n(1-|\lambda_n|)\gamma_n <\infty$ and $w\in L^1(\T)$
since $p>1$.
Associate with $\Lambda_0$
a second Carleson sequence $\Lambda_1=\{\lambda^\prime_n\}_n$
such that  the pseudo-hyperbolic distance between corresponding points
satisfies $|b_{\lambda^\prime_n}
(\lambda_n)|=e^{-\gamma_n^{1/p}}$.  Since
$\gamma_n\to\infty$ the elements of the sequence 
$\Lambda=\Lambda_0\cup\Lambda_1$  are
arbitrarily close and $\Lambda$ cannot be a Carleson sequence. By
construction,
$\log |B_\lambda(\lambda)|^{-1}\le P[w](\lambda)$ (as before, we
may possibly have to multiply
$u$ with some constant $c$ to have that condition also in the points
$\lambda^\prime_n$, but this operation conserves the integrability
condition), and
therefore $\Lambda\in\Int \hf$.
\end{example}

\providecommand{\bysame}{\leavevmode\hbox to3em{\hrulefill}\thinspace}
\providecommand{\MR}{\relax\ifhmode\unskip\space\fi MR }


\begin{thebibliography}{BRSHZE}

\bibitem[EiEs]{EiEs}
{\it Eiderman, V. Ya. \& Ess\'en, M.}
Uniqueness theorems for analytic and subharmonic functions (Russian),
Algebra i Analiz
\textbf{14} (2002), no. 6, 1--88 ;
    English translation in St.~Petersburg Math. J.
\textbf{14} (2002), no. 6.

\bibitem[Gar77]{Gar77}
{\it J. B. Garnett},
{Two remarks on interpolation by bounded analytic functions},
Banach spaces of analytic functions (Proc. Pelczynski Conf., Kent State
Univ., Kent, Ohio, 1976), pp. 32--40. Lecture Notes in Math., Vol. 604,
Springer, Berlin, 1977.

\bibitem[Gar81]{Gar}
{\it J.B. Garnett},
{Bounded analytic functions}, Academic Press, New York, 1981.

\bibitem[Har99]{Har}
{\it A. Hartmann},
{Free interpolation in Hardy-Orlicz spaces},
Studia Math. \textbf{135} (1999), no. 2, 179--190.

\bibitem[HaMa01]{HM2}
{\it A. Hartmann \& X. Massaneda},
{Interpolating sequences for holomorphic
functions of restricted growth}, Ill. J. Math. \textbf{46}
(2002), no.3, 929--945.

\bibitem[Hay]{Hay}
{\it W.K. Hayman},
{Identity theorems for functions of bounded characteristic}, J. London
Math.
Soc. \textbf{58} (1998), no. 2, 127--140.

\bibitem[He69]{He69}
{\it Heins, M.},
{Hardy Classes on Riemann Surfaces}, Springer Verlag, Berlin, Heidelberg,
New York, Lectures Notes in Mathematics no. 98, 1969.

\bibitem[HULL]{HULL}
Hiriart-Urruty, J.-B., Lemar\'echal, C. {\em Convex Analysis
and Minimization Algorithms I}, Grundlehren der mathematischen
Wissenschaften 305, Springer-Verlag, Berlin Heidelberg New York,
1993.

\bibitem[Lu96]{Lu96}
{\it D. Luecking},
{Zero sequences for Bergman spaces},
Complex Variables Theory Appl. \textbf{30} (1996), no. 4, 345--362.

\bibitem[LySe97]{LySe97}
{\it Lyubarskii, Yu. I. \& Seip, K.}
A uniqueness theorem for bounded analytic functions,
Bull. London Math. Soc.
\textbf{29} (1997), 49--52.

\bibitem[McC92]{McC92}
{\it J. McCarthy},
{Topologies on the Smirnov class},
J. Funct. Anal. \textbf{104} (1992), no. 1, 229--241.

\bibitem[Na56]{Na56}
{\it A.G. Naftalevi\v c},
{On interpolation by functions of
bounded characteristic (Russian)},  Vilniaus Valst. Univ. Moksl\c u Darbai.
Mat. Fiz. Chem. Moksl\c u \textbf{Ser. 5} (1956), 5--27.

\bibitem[NPT]{NPT}
{\it Nicolau, A., Pau, J., \& Thomas, P. J.}
Smallness sets for bounded holomorphic functions,
Journal d'Analyse Math\'ematique \textbf{82} (2000), 119--148.

\bibitem[NiTh]{NiTh}
{\it Nicolau, A., \& Thomas, P. J.}
Superharmonic envelopes for positive data on the disk, preprint, math/0304482.

\bibitem[Nik86]{niktr}
{\it N.K. Nikolski [Nikol'ski\u{\i}]},
{Treatise on the shift operator}, Springer-Verlag, Berlin etc., 1986.

\bibitem[Nik02]{Nik02}
\bysame
\emph{Operators, functions, and systems: an easy reading. Vol. 1,
Hardy, Hankel, and Toeplitz; Vol.2, Model Operators and Systems},
Mathematical Surveys and Monographs, 92 and 93.
American Mathematical Society, Providence, RI, 2002.


\bibitem[PaTh]{PaTh}
{\it Pau, J. \& Thomas, P. J.}
Decrease of bounded holomorphic
functions along discrete sets,
electronically available from math
arXiv :  http://arXiv.org/abs/math.CV/0106048, math/0106048

\bibitem[RosRov]{RR}
{\it M. Rosenblum \& J. Rovnyak},
\emph{Hardy classes and operator theory},
Oxford Mathematical Monographs. Oxford Science Publications.
The Clarendon Press, Oxford University
Press, New York, 1985.

\bibitem[Se93]{Se93}
{\it K. Seip}, {Beurling type density theorems in the unit disk},
Invent. Math., \textbf{113} (1993), 21--39.


\bibitem[ShHSh]{ShHSh}
{\it H.S. Shapiro \& A.L. Shields},
{On some interpolation problems for
analytic functions}, Amer. J. Math., \textbf{83} (1961), 513--532.

\bibitem[ShSh]{ShSh}
{\it J. Shapiro \& A. Shields},
{Unusual topological
properties of the Nevanlinna class}, Amer. J. Math. \textbf{97} (1975),
915--936.

\bibitem[Ya74]{yana2}
{\it N. Yanagihara},
{Interpolation theorems for the
class $N^+$}, Illinois J. Math., \textbf{18} (1974), 427--435.
\end{thebibliography}
\end{document}